\ifdefined\pdfminorversion
\fi
\documentclass[11pt]{amsart}
\usepackage[T1]{fontenc}
\usepackage{lmodern,microtype}
\usepackage{amsmath,amssymb,mathtools,amsthm}
\usepackage{geometry}
\usepackage{enumitem,booktabs,array,float,needspace}
\usepackage{tikz,tikz-cd}
\usetikzlibrary{arrows.meta,positioning,calc,decorations.markings}
\tikzset{graph vertex/.style={circle,fill=black,inner sep=1.5pt},%
         small graph vertex/.style={circle,fill=black,inner sep=1.2pt},%
         highlighted edge/.style={->,line width=1.15pt}}
\usepackage{aliascnt}
\usepackage[hidelinks]{hyperref}
\usepackage[nameinlink,noabbrev]{cleveref}
\numberwithin{equation}{section}
\newtheorem{theorem}{Theorem}[section]
\newtheorem{mainthm}{Theorem}

\crefname{mainthm}{Theorem}{Theorems}\Crefname{mainthm}{Theorem}{Theorems}
\newaliascnt{lemma}{theorem}\newtheorem{lemma}[lemma]{Lemma}\aliascntresetthe{lemma}
\newaliascnt{proposition}{theorem}\newtheorem{proposition}[proposition]{Proposition}\aliascntresetthe{proposition}
\newaliascnt{corollary}{theorem}\newtheorem{corollary}[corollary]{Corollary}\aliascntresetthe{corollary}
\theoremstyle{definition}
\newaliascnt{definition}{theorem}\newtheorem{definition}[definition]{Definition}\aliascntresetthe{definition}
\newaliascnt{example}{theorem}\newtheorem{example}[example]{Example}\aliascntresetthe{example}
\theoremstyle{remark}
\newaliascnt{remark}{theorem}\newtheorem{remark}[remark]{Remark}\aliascntresetthe{remark}
\crefname{theorem}{Theorem}{Theorems}\Crefname{theorem}{Theorem}{Theorems}
\crefname{lemma}{Lemma}{Lemmas}\Crefname{lemma}{Lemma}{Lemmas}
\crefname{proposition}{Proposition}{Propositions}\Crefname{proposition}{Proposition}{Propositions}
\crefname{corollary}{Corollary}{Corollaries}\Crefname{corollary}{Corollary}{Corollaries}
\crefname{definition}{Definition}{Definitions}\Crefname{definition}{Definition}{Definitions}
\crefname{remark}{Remark}{Remarks}\Crefname{remark}{Remark}{Remarks}
\crefname{example}{Example}{Examples}\Crefname{example}{Example}{Examples}
\newcommand{\Z}{\mathbb Z}

\newcommand{\id}{\operatorname{id}}
\newcommand{\Hom}{\operatorname{Hom}}
\newcommand{\toeq}{\mathrel{\overset{\scriptscriptstyle =}{\longrightarrow}}}
\newcommand{\lefteq}{\mathrel{\overset{\scriptscriptstyle =}{\longleftarrow}}}
\newcommand{\pC}{\pi^C}
\newcommand{\Lt}{L_{\mathrm{tail}}}
\newcommand{\Sn}{S^k_{2n_1\times\cdots\times2n_k}}
\newcommand{\ev}{\operatorname{ev}}
\newcommand{\nn}{\mathbf n}
\setlist[enumerate]{label=\textup{(\roman*)},leftmargin=2.1em,itemsep=3pt}
\title[Sphere homotopy groups and the GLMY loop digraph]{Sphere Homotopy Groups of Digraphs and the Vanishing of the Original Higher GLMY Homotopy Groups}
\author[Jingyan Li]{Jingyan Li$^{*}$}
\address[Jingyan Li]{Beijing Institute of Mathematical Sciences and Applications\newline
No.~544 Hefangkou Village, Huairou, Beijing 101408, China}
\email{jingyanli@bimsa.cn}
\thanks{$^{*}$Corresponding author. E-mail: \href{mailto:jingyanli@bimsa.cn}{jingyanli@bimsa.cn}.}
\author{Jie Wu}
\address[Jie Wu]{Beijing Institute of Mathematical Sciences and Applications\newline
No.~544 Hefangkou Village, Huairou, Beijing 101408, China}
\email{wujie@bimsa.cn}
\author{Shing-Tung Yau}
\address[Shing-Tung Yau]{Yau Mathematical Sciences Center\newline
Jing Zhai, Tsinghua University, Haidian District, Beijing 100084, China}
\email{styau@tsinghua.edu.cn}
\makeatletter
\renewcommand{\@setaddresses}{%
  \par\begingroup\normalfont\normalsize
  \def\author##1{}%
  \def\address##1##2{%
    \par\addvspace{\bigskipamount}\noindent
    ##1\par\nobreak\noindent##2\par}%
  \def\email##1##2{%
    \nobreak\noindent E-mail: \href{mailto:##2}{\textcolor{blue}{##2}}\par}%
  \addresses
  \endgroup}
\makeatother
\date{}
\hypersetup{pdftitle={Sphere Homotopy Groups of Digraphs and the Vanishing of the Original Higher GLMY Homotopy Groups},pdfauthor={Jingyan Li, Jie Wu, and Shing-Tung Yau}}
\begin{document}
\raggedbottom
\begin{abstract}
We prove that $\pi_k^{\mathrm{GLMY}}(G,*)=0$ for every based digraph
$(G,*)$ and every $k\ge2$, where the groups are defined using the
original loop digraph of Grigor'yan, Lin, Muranov and Yau.
We define sphere homotopy groups $\pC_k(G,*)$ using based maps
from digraph spheres presented as boundary quotients of boxes of
oriented lines. We require these quotients to have the integral GLMY
path homology of a sphere. In degree one, the sphere homotopy group
is naturally isomorphic to the GLMY fundamental group. For every
$k\ge1$, we construct a natural Hurewicz homomorphism
$h_k^C:\pC_k(G,*)\to H_k(G;\Z)$ to integral GLMY path homology.
For $k\ge2$, the sphere homotopy groups of a standard digraph $k$-sphere
agree with the homotopy groups of the topological $k$-sphere in every
positive degree. In degree $k$, the Hurewicz homomorphism is an
isomorphism and sends the class represented by the standard quotient
map to the oriented generator. The sphere homotopy groups of the
six-vertex inverse suspension of the alternating four-cycle agree
with the homotopy groups of the topological $2$-sphere in every
positive degree.
We recover the sphere homotopy groups by iterating a loop digraph
whose vertices are based loops on a common infinite alternating
line. Each loop equals the basepoint for all sufficiently large
parameter values.
\end{abstract}
\maketitle

\section{Introduction}\label{sec:intro}

We study two constructions of homotopy groups for based digraphs.
The first iterates the original GLMY loop digraph, and the second
uses based maps from digraph spheres.

For a based digraph $(G,*)$, Grigor'yan, Lin, Muranov and Yau
(GLMY) defined the fundamental group $\pi_1(G,*)$ using the
loop digraph $L(G,*)$. Its vertices are parametrized based loops.
The components of its underlying undirected graph form
$\pi_1(G,*)$ under concatenation. They then defined $\pi_2^{\mathrm{GLMY}}(G,*)$ to be the
fundamental group of $L(G,*)$ at the length-zero loop $e_{G,*}$,
and defined the higher groups by iteration
\cite[Definitions~4.8, 4.10 and 4.24]{GLMY}.
In \cite[Remark~4.27]{GLMY}, the authors observed that changing the
arrows between loops can preserve the components, while the
fundamental group of the loop digraph may change.

\begin{mainthm}\label{thm:GLMY-higher-vanishing}
For every based digraph $(G,*)$, the original GLMY loop
digraph satisfies
\[
 \pi_1\bigl(L(G,*),e_{G,*}\bigr)=0.
\]
Consequently, $\pi_k^{\mathrm{GLMY}}(G,*)=0$ for every $k\ge2$.
\end{mainthm}

The descriptions of topological homotopy groups by sphere maps
and iterated loop spaces
\cite[Sections~4.1 and~4.3]{Hatcher} motivate the following
construction for digraphs.
We collapse the full coordinate boundary of a box product of $k$ finite
oriented lines to a point. Boundary quotients of boxes with identical
factors are the digraph spheres of
\cite[Definition~1.3.9 and Corollary~1.3.10]{EldridgeEtAl}.
We allow the oriented line factors to vary and call the quotient a
$k$-dimensional digraph sphere if its integral GLMY path homology is
isomorphic to the singular homology $H_*(S^k;\Z)$.

For a based digraph $(G,*)$, we define
its $k$th sphere homotopy group $\pC_k(G,*)$ as the set of $C$-homotopy
classes of representatives $f=\bar f\circ q$. Here $q$ is a chosen
boundary quotient map onto a $k$-dimensional digraph sphere, and
$\bar f$ is a based digraph map from that sphere to $(G,*)$.
In degree one this construction recovers the GLMY fundamental group
$\pi_1(G,*)$ (\cref{prop:pi1}).
These groups are naturally isomorphic to the cubical homotopy groups
$A_k(G)$ of Eldridge et al.\ (\cref{rem:cubical-comparison}).

For $\nn=(n_1,\ldots,n_k)$ with $n_i\ge2$, let
$B_{\nn}=J_{2n_1}\Box\cdots\Box J_{2n_k}$ be the standard alternating box
and let $q_{\nn}:(B_{\nn},\partial B_{\nn})\to(\Sn,*)$ be its
boundary quotient map (Section~\ref{sec:spheres}).
For a standard representative $f=\bar f\circ q_{\nn}$, the Hurewicz
homomorphism $h_k^C:\pC_k(G,*)\to H_k(G;\Z)$ sends $[f]$ to the image
of the oriented top class of $\Sn$ under $H_k(\bar f)$
(\cref{thm:top-map}).

Eldridge et al. compare the cubical homotopy groups of a
digraph with the homotopy groups of its cubical nerve realization
\cite[Theorem~2.4.5 and Proposition~2.5.3]{EldridgeEtAl}.
Kishimoto and Tong establish a realization theorem for the discrete
homotopy groups of quasisymmetric cubical sets
\cite[Theorem~1.1]{KishimotoTong}. We identify the $C$-homotopy classes
of sphere representatives with spherical grid classes in the cubical
nerve and show that this correspondence preserves multiplication.
We then apply the realization theorem
(\cref{prop:sphere-grid-realization}) to identify the sphere homotopy
groups of standard digraph spheres and of the six-vertex inverse
suspension with the homotopy groups of the corresponding topological
spheres.

\begin{mainthm}\label{cor:standard-sphere-separation}
Let $k\ge2$ and $n_1,\ldots,n_k\ge2$, and let $s_0$ be a basepoint
of the topological sphere $S^k$. For every $j\ge1$,
\[
 \pC_j(\Sn,*)\cong\pi_j(S^k,s_0).
\]
In degree $k$, the Hurewicz homomorphism
\[
 h_k^C:\pC_k(\Sn,*)\longrightarrow H_k(\Sn;\Z)
\]
is an isomorphism and sends $[q_{\nn}]$ to the oriented generator.
\end{mainthm}

Taking $n_1=\cdots=n_k=2$ gives a standard sphere with $3^k+1$
vertices in each dimension $k\ge2$. Its original GLMY homotopy
group in the same degree is trivial by \cref{thm:GLMY-higher-vanishing}.

Grigor'yan, Muranov and Yau computed the GLMY path homology of
directed suspensions and exhibited a top cycle for the octahedral
suspension, with arrows from the rim to the poles
\cite[Proposition~3.7 and Example~3.8]{GMY}.
Let $(G,a)$ be the six-vertex inverse suspension of the alternating
four-cycle, with arrows from both poles to every rim vertex.
Let $F:(J_4\Box J_4,\partial(J_4\Box J_4))\to(G,a)$ be the explicit
map of pairs in \cref{prop:octahedral-sphere-class}.
For any basepoint $s_0$ of $S^2$, \cref{cor:six-vertex-separation} gives
\[
 \pC_j(G,a)\cong\pi_j(S^2,s_0)\quad(j\ge1),\qquad
 \pC_2(G,a)=\Z[F].
\]
In degree two, $h_2^C$ is an isomorphism that sends $[F]$ to the
generator $[w]$ of $H_2(G;\Z)$, where $w$ is the cycle in
\eqref{eq:rim-chain}.

The six-vertex example rules out a recursion for sphere homotopy
groups using the original, replenished, or simple loop digraph
(Section~\ref{sec:tail}).
We use the pointwise-tail loop digraph $\Lt(G,*)$, a model of the
Eldridge et al. loop digraph for the right-shrinking sequence
\cite[Definition~1.4.1]{EldridgeEtAl}. Slicing in the last coordinate
of a standard box and evaluation give natural isomorphisms
(\cref{thm:tail-recursion})
\[
 \pC_k(G,*)\cong\pC_{k-1}\bigl(\Lt(G,*)\bigr)
 \qquad(k\ge2).
\]

Grid maps also occur in $A$-homotopy theory for simplicial
complexes and graphs \cite{BBLL,BKLW,BarceloLaubenbacher}.
When every arrow has its reverse, the sphere homotopy groups recover the
$A$-homotopy groups of the underlying simple graph
(\cref{rem:symmetric-targets}). For simple graphs, Lutz
detected higher $A$-homotopy classes by a map to discrete cubical
singular homology \cite{Lutz}, and Carranza and Kapulkin gave a
cubical model and a Hurewicz theorem \cite{CarranzaKapulkin}.
The homology used here is integral GLMY path homology of the
digraph. Hom-complex approaches to graph and digraph homotopy
are developed in \cite{DochtermannEJC,DochtermannJCTA,DochtermannSingh}.

Section~\ref{sec:preliminaries} recalls the GLMY definitions, and
Section~\ref{sec:vanishing} proves \cref{thm:GLMY-higher-vanishing}.
Sections~\ref{sec:spheres}--\ref{sec:groups} construct the sphere
homotopy groups. Section~\ref{sec:top-class} constructs the Hurewicz
homomorphism. Section~\ref{sec:examples} proves
\cref{cor:standard-sphere-separation} and gives the six-vertex example.
Section~\ref{sec:tail} proves the loop recursion.
Appendix~\ref{app:homology} computes the integral GLMY path homology
of the standard box quotients.

\Needspace{9\baselineskip}
\section{Preliminaries}\label{sec:preliminaries}

We recall digraph homotopy and the original GLMY loop digraph.
The proof in Section~\ref{sec:vanishing} uses the shrinking rule
for arrows between based loops.

\subsection{Based digraphs, shrinking maps and homotopy}

\begin{definition}\label{def:digraph-map}
A \emph{digraph} is a pair $G=(V_G,E_G)$, where
\[
 E_G\subseteq\{(x,y)\in V_G\times V_G:x\ne y\}.
\]
We write $x\to y$ for $(x,y)\in E_G$ and
\[
 x\toeq y
 \quad\Longleftrightarrow\quad
 x=y\text{ or }x\to y.
\]
We also write $x\lefteq y$ for $y\toeq x$.
For a digraph $H$, a vertex map $f:V_G\to V_H$ is a \emph{digraph map} if
\[
 x\to y\text{ in }G\quad\Longrightarrow\quad f(x)\toeq f(y)\text{ in }H.
\]
Thus a digraph map may collapse an arrow.
A \emph{digraph homomorphism} is a digraph map that collapses no arrow.
\end{definition}

\begin{definition}\label{def:based-pair}
A \emph{based digraph} is a pair $(G,*)$ consisting of a digraph $G$ and a
distinguished vertex $*\in V_G$.  A \emph{based digraph map}
$f:(G,*)\to(H,*)$ between based digraphs is a digraph map $f:G\to H$
satisfying $f(*)=*$.

Let $X$ be a digraph and $A$ a subdigraph of $X$.  A digraph map of pairs
\[
 f:(X,A)\longrightarrow(G,*)
\]
is a digraph map $f:X\to G$ satisfying $f(a)=*$ for every $a\in V_A$.

We write $\Hom((X,A),(G,*))$ for the set of these digraph maps of pairs.
All homotopy groups in this paper are based, and homotopies of digraph maps of pairs
are relative to the displayed subdigraph unless stated otherwise.
\end{definition}

\begin{definition}\label{def:oriented-line}
For $m\ge0$, a \emph{finite oriented line} $I_m$ has vertex set
$\{0,1,\ldots,m\}$ and, for each $0\le i<m$, exactly one of the two arrows
\[
 i\to i+1,
 \qquad
 i+1\to i.
\]
There are no other arrows in $I_m$.
Its boundary $\partial I_m$ is the induced subdigraph on $\{0,m\}$,
with one vertex when $m=0$. For a based digraph $(G,*)$, a based loop on
$I_m$ is a digraph map of pairs
\[
 \gamma:(I_m,\partial I_m)\longrightarrow(G,*).
\]
\end{definition}

\begin{definition}\label{def:line-shrinking}
Let $I_N$ and $I_m$ be finite oriented lines.  A \emph{shrinking}
\[
 h:I_N\longrightarrow I_m
\]
is a surjective nondecreasing digraph map satisfying
\[
 h(0)=0,
 \qquad
 h(N)=m.
\]
Necessarily $N\ge m$.
\end{definition}

\begin{definition}\label{def:box-product}
For digraphs $X$ and $Y$, the \emph{box product} $X\Box Y$ has vertex set
$V_X\times V_Y$ and arrows
\[
 (x,y)\to(x',y)\quad\text{if }x\to x',
 \qquad
 (x,y)\to(x,y')\quad\text{if }y\to y'.
\]
Thus every arrow changes exactly one coordinate.
\end{definition}

\begin{definition}\label{def:pointwise-comparison}
Let $f_0,f_1:X\to G$ be digraph maps with the same source.  We write
\[
 f_0\toeq f_1
\]
if
\[
 f_0(x)\toeq f_1(x)
 \qquad\text{for every }x\in V_X.
\]
We write $f_0\lefteq f_1$ for $f_1\toeq f_0$.
\end{definition}

\begin{definition}\label{def:homotopy}
Digraph maps $f_0,f_1:X\to G$ are \emph{homotopic} if there exist a
finite oriented line $I_m$ and a digraph map
\[
 H:X\Box I_m\longrightarrow G
\]
such that $H(x,0)=f_0(x)$ and $H(x,m)=f_1(x)$ for every $x\in V_X$.
For $m=1$, this is a \emph{one-step direct homotopy} when
$I_1=[0\to1]$, and a \emph{one-step inverse homotopy} when
$I_1=[0\leftarrow1]$.  The respective conditions are
$f_0\toeq f_1$ and $f_1\toeq f_0$.
Equivalently, a homotopy is a finite sequence of one-step direct or
inverse homotopies.
We write $f_0\simeq f_1$ for homotopy.

For digraph maps of pairs $f_0,f_1:(X,A)\to(G,*)$, a homotopy $H$ from $f_0$ to $f_1$
is \emph{relative to $A$} if
\[
 H:(X\Box I_m,A\Box I_m)\longrightarrow(G,*)
\]
is a digraph map of pairs.
\end{definition}

\subsection{The GLMY loop digraph and homotopy groups}

We use the original GLMY definitions
\cite[Definitions~4.4, 4.5, 4.8, 4.10 and 4.24]{GLMY}.
We allow infinite targets. The images of finitely many based loops
lie in a finite induced subdigraph of $G$ containing $*$.

\begin{definition}\label{def:classical-loop}
For a based digraph $(G,*)$, the original GLMY loop digraph $L(G,*)$
has as vertices all based loops
\[
 \gamma:(I_m,\partial I_m)\longrightarrow(G,*)
\]
where $I_m$ ranges over finite oriented lines. The constant loop
\mbox{$e_{G,*}:(I_0,\partial I_0)\to(G,*)$} is the basepoint.
For distinct based loops $\gamma:(I_m,\partial I_m)\to(G,*)$ and
$\delta:(I_n,\partial I_n)\to(G,*)$, an arrow $\gamma\to\delta$ exists
if and only if a shrinking $h$ satisfies
\begin{equation}\label{eq:classical-arrow}
\begin{aligned}
 &m\ge n:\quad h:I_m\to I_n,
 &&\gamma(t)\toeq\delta(h(t))\quad(0\le t\le m);\\
 &m<n:\quad h:I_n\to I_m,
 &&\gamma(h(t))\toeq\delta(t)\quad(0\le t\le n).
\end{aligned}
\end{equation}
\end{definition}

\begin{definition}\label{def:pi0}
For a digraph $H$, let $\pi_0(H)$ denote the set of connected components of
the underlying undirected graph. Equivalently, $x,y\in V_H$ lie in the
same component if there is a finite oriented line $I_m$ and a
digraph map $I_m\to H$ sending $0$ to $x$ and $m$ to $y$.
We call $H$ \emph{weakly connected} if $\pi_0(H)$ is a singleton.
\end{definition}

Two based loops are $C$-homotopic if they lie in the same component of
$L(G,*)$. The GLMY fundamental group is $\pi_1(G,*)=\pi_0(L(G,*))$, with
multiplication induced by loop concatenation. For paths $\gamma:I_m\to G$
and $\delta:I_n\to G$ with $\gamma(m)=\delta(0)$, concatenate their
parameter lines by identifying $m$ in the first with $0$ in the second
and relabeling vertex $t$ in the second by $m+t$. Each arrow keeps its
orientation. On this line of length $m+n$, set
\[
 (\gamma\vee\delta)(t)=
 \begin{cases}
  \gamma(t),&0\le t\le m,\\
  \delta(t-m),&m\le t\le m+n.
 \end{cases}
\]
For based loops this is the GLMY concatenation
\cite[Definition~4.16(ii)]{GLMY}. The length-zero loop is a strict unit:
\begin{equation}\label{eq:loop-strict-unit}
 e_{G,*}\vee\gamma=\gamma=\gamma\vee e_{G,*}
 \qquad\text{for every based loop }\gamma.
\end{equation}
These equalities include the oriented parameter lines. We take the empty
concatenation to be $e_{G,*}$.
By \eqref{eq:classical-arrow}, for a based loop $\gamma$ and a
shrinking $h$, the loops $\gamma$ and $\gamma\circ h$ represent
the same class in $\pi_1(G,*)$.

\begin{definition}\label{def:GLMY-recursion}
Let $(G,*)$ be a based digraph.
The iterated based loop digraphs are defined recursively by
\[
 L^0(G,*)=(G,*),\qquad
 L^{r+1}(G,*)=L\bigl(L^r(G,*)\bigr)
 \quad(r\ge0).
\]
For each integer $r\ge0$, the basepoint of $L^{r+1}(G,*)$ is the constant
loop on $I_0$ at the basepoint of $L^r(G,*)$.  For $k\ge1$, define
\begin{equation}\label{eq:GLMY-recursive}
 \pi_k^{\mathrm{GLMY}}(G,*)
 =\pi_1\bigl(L^{k-1}(G,*)\bigr).
\end{equation}
We omit the superscript $\mathrm{GLMY}$ in degree one.
In particular,
\[
 \pi_2^{\mathrm{GLMY}}(G,*)=\pi_1\bigl(L(G,*)\bigr).
\]
\end{definition}

Postcomposition with a based digraph map $f:(G,*)\to(H,*)$ preserves the pointwise comparisons in
\eqref{eq:classical-arrow} for the same shrinkings. It therefore induces
a based digraph map
\[
 L(f):L(G,*)\longrightarrow L(H,*),\qquad \gamma\longmapsto f\circ\gamma.
\]
Thus $L$ is a functor on based digraphs, and so are its iterates.
 
\section{Vanishing of the original higher GLMY homotopy groups}\label{sec:vanishing}

We compare two staircases in a finite box of oriented lines.
Their images under loop concatenation relate a given loop in
$L(G,*)$ to an auxiliary loop that can be contracted.

\begin{proof}[Proof of \cref{thm:GLMY-higher-vanishing}]
We first prove that every based loop in $\bigl(L(G,*),e_{G,*}\bigr)$
is $C$-null-homotopic. Let
\[
 \Gamma:(I_m,\partial I_m)\longrightarrow
       \bigl(L(G,*),e_{G,*}\bigr)
\]
be a based loop. If $m=0$, the loop is already constant.
Assume $m\ge1$, and write
\[
 \gamma_i=\Gamma(i),\qquad
 \gamma_i:(I_{\ell_i},\partial I_{\ell_i})\longrightarrow(G,*)
 \quad(0\le i\le m).
\]
The line $I_m$ parametrizes $\Gamma$. Each value $\gamma_i$ is a loop
in $G$ with its own oriented parameter line $I_{\ell_i}$. At the two endpoints,
$\gamma_0=\gamma_m=e_{G,*}$ and $\ell_0=\ell_m=0$.

\medskip
Fix $1\le i\le m$. We choose a finite oriented line $K_i$, shrinkings
\[
 \sigma_i:K_i\longrightarrow I_{\ell_{i-1}},
 \qquad
 \tau_i:K_i\longrightarrow I_{\ell_i},
\]
and a based loop $\beta_i:(K_i,\partial K_i)\to(G,*)$ such that
\begin{equation}\label{eq:GLMY-block-upper}
 \gamma_{i-1}\circ\sigma_i\toeq\beta_i,
 \qquad
 \gamma_i\circ\tau_i\toeq\beta_i.
\end{equation}
If $\gamma_{i-1}=\gamma_i$, take $K_i=I_{\ell_i}$,
$\sigma_i=\tau_i=\id$, and $\beta_i=\gamma_i$.
If the loops are distinct, the arrow between $i-1$ and $i$ in $I_m$
gives an arrow between the corresponding loops in $L(G,*)$, in the
same direction. Take $K_i$ to be the longer parameter line and use
the identity on this line and the shrinking in
\eqref{eq:classical-arrow} to the shorter line. At equal lengths
the shrinking is the identity on the vertex labels. It preserves
each arrow, so the two oriented lines agree. Thus, if $i-1\to i$,
the two maps on $K_i$ satisfy
$\gamma_{i-1}\circ\sigma_i\toeq\gamma_i\circ\tau_i$, and we take
$\beta_i=\gamma_i\circ\tau_i$. If $i\to i-1$, the comparison is
$\gamma_i\circ\tau_i\toeq\gamma_{i-1}\circ\sigma_i$, and we take
$\beta_i=\gamma_{i-1}\circ\sigma_i$. In either case, one comparison
in \eqref{eq:GLMY-block-upper} is equality and the other is the
comparison given by the arrow.

Let $\beta=\beta_1\vee\cdots\vee\beta_m$ be the based loop in $(G,*)$
on the line obtained by concatenating $K_1,\ldots,K_m$. For every choice
$\delta_i\in\{\gamma_{i-1},\gamma_i\}$ we have
\begin{equation}\label{eq:GLMY-concatenated-upper}
 \delta_1\vee\cdots\vee\delta_m\toeq\beta
 \qquad\text{in }L(G,*).
\end{equation}
For each $i$, take $h_i=\sigma_i$ when $\delta_i=\gamma_{i-1}$ and
$h_i=\tau_i$ when $\delta_i=\gamma_i$. The two choices agree when
$\gamma_{i-1}=\gamma_i$. Write $N_i$ for the
length of $K_i$ and $M_i$ for the length of the line of $\delta_i$.
The formula
\begin{equation}\label{eq:GLMY-concatenated-shrinking}
 h\left(\sum_{j<i}N_j+t\right)
 =\sum_{j<i}M_j+h_i(t),
 \qquad 0\le t\le N_i,\quad 1\le i\le m,
\end{equation}
defines a vertex map from the parameter line of $\beta$ to that of
$\delta_1\vee\cdots\vee\delta_m$. At the junction of blocks $i$ and
$i+1$, the two formulas give the same value:
\[
 \sum_{j<i} M_j+h_i(N_i)=\sum_{j\le i} M_j
 =\sum_{j<i+1} M_j+h_{i+1}(0).
\]
If $N_i=0$, then $M_i=0$, so any coincident junctions give the same
value as well. On each block, $h$ is nondecreasing and maps onto the
corresponding target block. The target blocks occur in order and
cover the target line, so $h$ is a nondecreasing surjection. Every
arrow of the source line belongs to one of its blocks and is
preserved or collapsed by that block's $h_i$. Hence $h$ is a
shrinking. At the source vertex with local parameter $t$ in block
$i$, the two loop values to compare are $\delta_i(h_i(t))$ and
$\beta_i(t)$. They satisfy $\delta_i(h_i(t))\toeq\beta_i(t)$ by
\eqref{eq:GLMY-block-upper}. Since $\sum_i N_i\ge\sum_i M_i$,
this is exactly the comparison in \eqref{eq:classical-arrow} for
\eqref{eq:GLMY-concatenated-upper}, with equality allowed.

\medskip
For each vertex $(t_1,\ldots,t_m)$ of $I_m^{\Box m}$, concatenate the
loops $\Gamma(t_j)$ in increasing order of the coordinate index $j$,
using each loop's own parameter line. The resulting map is
\[
 \Phi:I_m^{\Box m}\longrightarrow L(G,*),\qquad
 \Phi(t_1,\ldots,t_m)
 =\gamma_{t_1}\vee\cdots\vee\gamma_{t_m}.
\]
Suppose a box arrow changes the $j$th coordinate from $t_j$ to $t_j'$.
Then $t_j\to t_j'$ in $I_m$, and $\Gamma$ gives
$\gamma_{t_j}\toeq\gamma_{t_j'}$. If these loops agree, the
concatenations agree. Otherwise, use the shrinking that witnesses
this arrow on block $j$ and the identity maps on all other blocks.
The total lengths of the concatenations differ by the difference
between $\ell_{t_j}$ and $\ell_{t_j'}$. The block formula
\eqref{eq:GLMY-concatenated-shrinking} therefore gives a shrinking
from the longer concatenated line to the shorter one, with the
identity on vertex labels in the equal-length case. On the
unchanged blocks the loop values agree after this alignment. On the
changed block they have the comparison in the direction
$\gamma_{t_j}\to\gamma_{t_j'}$. Thus \eqref{eq:classical-arrow}
gives the comparison from the image of the source vertex to that
of the target vertex. This proves that $\Phi$ is a digraph map.

Start at $(0,1,\ldots,m-1)$. At step $i$, for $i=1,\ldots,m$, change
coordinate $i$ from $i-1$ to $i$ and keep the others fixed. The final
vertex is $(1,2,\ldots,m)$. The box arrow between the vertices at steps
$i-1$ and $i$ has the direction of the edge between $i-1$ and $i$ in
the $i$th copy of $I_m$. Thus the vertex sequence defines a digraph
map on the original line $I_m$, including when that edge points from
$i$ to $i-1$. Let $S:I_m\to L(G,*)$ be its image under $\Phi$:
\begin{equation}\label{eq:GLMY-staircase-loop}
 S(i)
 =\Phi\bigl(\underbrace{1,\ldots,i}_{i\ \text{coordinates}},
            \underbrace{i,\ldots,m-1}_{m-i\ \text{coordinates}}\bigr),
 \qquad 0\le i\le m.
\end{equation}
At step $i$ the coordinates numbered $1$ through $i$ have already
increased, and the rest still have their initial values. The first
list is empty at $i=0$, and the second at $i=m$. Since $\gamma_0$ and
$\gamma_m$ are both the strict unit in \eqref{eq:loop-strict-unit},
\[
 \begin{aligned}
 S(0)&=e_{G,*}\vee\gamma_1\vee\cdots\vee\gamma_{m-1},\\
 S(m)&=\gamma_1\vee\cdots\vee\gamma_{m-1}\vee e_{G,*}.
 \end{aligned}
\]
Both values are $\rho:=\gamma_1\vee\cdots\vee\gamma_{m-1}$, which
equals $e_{G,*}$ for $m=1$. Thus $S$ is a based loop in
$(L(G,*),\rho)$. The $j$th factor of $S(i)$ is $\gamma_j$ for
$j\le i$ and $\gamma_{j-1}$ for $j>i$. By
\eqref{eq:GLMY-concatenated-upper},
\begin{equation}\label{eq:GLMY-staircase-upper}
 S(i)\toeq\beta\quad(0\le i\le m),
 \qquad \rho\toeq\beta.
\end{equation}

The common target $\beta$ need not equal $\rho$. To keep the
basepoint fixed while moving the interior of $S$ to $\beta$,
we first repeat the value $\rho$ at each end. Let
$\widetilde I_{m+2}$ have vertices $0,\ldots,m+2$.
Its subline on $1,\ldots,m+1$ is obtained from $I_m$ by adding
$1$ to each vertex label, and its two additional arrows are
$0\to1$ and $m+2\to m+1$.
The map
\[
 r:\widetilde I_{m+2}\longrightarrow I_m,\qquad
 r(t)=\min\{m,\max\{0,t-1\}\},
\]
preserves the translated arrows on the subline $1,\ldots,m+1$.
It collapses $0\to1$ to $0$ and $m+2\to m+1$ to $m$. Its vertex map
is nondecreasing and surjective, so $r$ is a shrinking. Consequently
$S\circ r$ represents the same $C$-homotopy class as $S$.
Using $J_2=[0\to1\leftarrow2]$, define
\[
 \mathcal H:
 \bigl(\widetilde I_{m+2}\Box J_2,
       \partial\widetilde I_{m+2}\Box J_2\bigr)
 \longrightarrow\bigl(L(G,*),\rho\bigr)
\]
by
\begin{equation}\label{eq:GLMY-staircase-contraction}
 \mathcal H(t,s)=
 \begin{cases}
  \rho,&t\in\{0,m+2\},\\
  S(r(t)),&1\le t\le m+1,\ s=0,\\
  \beta,&1\le t\le m+1,\ s=1,\\
  \rho,&1\le t\le m+1,\ s=2.
 \end{cases}
\end{equation}
At $s=0$ the slice is $S\circ r$, and at $s=2$ it is constant at
$\rho$. At $s=1$, each arrow between interior vertices collapses,
and the two endpoint arrows give the comparison $\rho\toeq\beta$.
For $1\le t\le m+1$, the time arrows $0\to1$ and $2\to1$ give,
respectively, $S(r(t))\toeq\beta$ and $\rho\toeq\beta$, which hold
by \eqref{eq:GLMY-staircase-upper}. At $t=0,m+2$ the time arrows
collapse to $\rho$. Therefore $\mathcal H$ is a digraph map of
pairs. Using identity shrinkings on $\widetilde I_{m+2}$, the slices
at $s=0,1,2$ give a one-step direct $C$-homotopy followed by a
one-step inverse $C$-homotopy from $S\circ r$ to the constant loop
at $\rho$. The constant loop on
$\widetilde I_{m+2}$ is the precomposition of the length-zero loop
at $\rho$ with the shrinking to $I_0$. Hence
\[
 [S]=1\quad\text{in }\pi_1\bigl(L(G,*),\rho\bigr).
\]

\medskip
We now construct two staircases from $(0,\ldots,0)$ to $(1,2,\ldots,m)$
in $I_m^{\Box m}$. A staircase increases one coordinate label by $1$
at each step. We orient each edge of its parameter line as the
corresponding edge of the box.

For the first staircase, start at $(0,\ldots,0)$ and, for
$j=2,\ldots,m$ in increasing order, increase coordinate $j$ from
$0$ to $j-1$. This segment ends at $(0,1,\ldots,m-1)$.
The segment has $\binom{m}{2}$ edges, whose orientations give its
parameter line $I^P_{\binom{m}{2}}$. Its image under $\Phi$ is the
digraph map
\[
 P:I^P_{\binom{m}{2}}\longrightarrow L(G,*),
 \qquad P(0)=e_{G,*},\qquad P\Bigl(\binom{m}{2}\Bigr)=\rho.
\]
Continue along the staircase defining $S$ in
\eqref{eq:GLMY-staircase-loop} to $(1,2,\ldots,m)$. Under $\Phi$,
the full first staircase maps to $P\vee S$.

For the second staircase, start again at $(0,\ldots,0)$ and increase
only the last coordinate from $0$ to $m$. This segment is the digraph
map $t\mapsto(0,\ldots,0,t)$ on $I_m$, and it ends at
$(0,\ldots,0,m)$. Since $\gamma_0=e_{G,*}$, \eqref{eq:loop-strict-unit} gives
\[
 \Phi(0,\ldots,0,t)
 =\underbrace{\gamma_0\vee\cdots\vee\gamma_0}_{m-1\ \text{factors}}
   \vee\gamma_t
 =\gamma_t=\Gamma(t)
 \qquad (0\le t\le m).
\]
Thus this segment has image $\Gamma$ on its original parameter line.
For the remaining segment, keep the last coordinate equal to $m$
and, for $j=1,\ldots,m-1$ in increasing order, increase coordinate
$j$ from $0$ to $j$. This is the image of the first segment of the
first staircase under
$(0,t_2,\ldots,t_m)\mapsto(t_2,\ldots,t_m,m)$. In the first segment,
a step changes some $t_j$ by one, for $2\le j\le m$. In the
translated segment the same label changes by one in coordinate
$j-1$. Both coordinates are copies of $I_m$, so the two steps have
the same arrow direction in their parameter lines. Moreover,
since $\gamma_0=\gamma_m=e_{G,*}$, \eqref{eq:loop-strict-unit} gives
\[
 \Phi(0,t_2,\ldots,t_m)
 =\gamma_{t_2}\vee\cdots\vee\gamma_{t_m}
 =\Phi(t_2,\ldots,t_m,m).
\]
The remaining segment therefore has the same parametrized image
$P$, including its oriented parameter line. The full second staircase
ends at $(1,2,\ldots,m)$ and maps to $\Gamma\vee P$. When $m=1$,
both segments with image $P$ have length zero and $\rho=e_{G,*}$.

\Cref{fig:GLMY-staircases} shows the two staircases for
$I_3=[0\to1\leftarrow2\to3]$.

\begin{figure}[H]
\centering
\begin{tikzpicture}[
  x={(18.45mm,-4.7mm)},y={(10.9mm,5.1mm)},z={(0mm,10.4mm)},
  scale=.82,
  lattice/.style={draw=black!31,line width=.25pt,shorten <=.65pt,shorten >=.65pt},
  lattice rim/.style={lattice,draw=black!46,line width=.32pt},
  lattice direction/.style={postaction={decorate},decoration={markings,
    mark=at position #1 with {\arrow{Stealth[length=2.45pt,width=2.05pt]}}}},
  route underlay/.style={draw=white,line width=2.25pt,line cap=round},
  route top/.style={line width=1.05pt,shorten <=.9pt,shorten >=1.7pt,
    -{Stealth[length=3.8pt,width=3.3pt]}},
  pathP/.style={draw=blue!67!black},
  pathS/.style={draw=teal!68!black},
  pathG/.style={draw=red!59!black},
  path label/.style={font=\small,fill=white,inner sep=1.2pt}
]
\foreach \a/\b in {0/1,2/1,2/3}{
  \foreach \i in {0,1,2,3}{\foreach \j in {0,1,2,3}{
    \ifnum\numexpr\i*(3-\i)+\j*(3-\j)\relax=0
      \def\latticeweight{lattice rim}
    \else
      \def\latticeweight{lattice}
    \fi
    \draw[\latticeweight,lattice direction=.76] (\a,\i,\j)--(\b,\i,\j);
    \draw[\latticeweight,lattice direction=.53] (\i,\a,\j)--(\i,\b,\j);
    \draw[\latticeweight,lattice direction=.55] (\i,\j,\a)--(\i,\j,\b);
  }}
}
\foreach \routepass in {route underlay,route top}{
  \foreach \source/\target/\routecolor in {
    {(0,0,0)}/{(0,1,0)}/pathP,
    {(0,1,0)}/{(0,1,1)}/pathP,
    {(0,1,2)}/{(0,1,1)}/pathP,
    {(0,1,2)}/{(1,1,2)}/pathS,
    {(1,2,2)}/{(1,1,2)}/pathS,
    {(1,2,2)}/{(1,2,3)}/pathS,
    {(0,0,0)}/{(0,0,1)}/pathG,
    {(0,0,2)}/{(0,0,1)}/pathG,
    {(0,0,2)}/{(0,0,3)}/pathG,
    {(0,0,3)}/{(1,0,3)}/pathP,
    {(1,0,3)}/{(1,1,3)}/pathP,
    {(1,2,3)}/{(1,1,3)}/pathP%
  }{
    \draw[\routecolor,\routepass] \source--\target;
  }
}
\foreach \i in {0,1,2,3}{\foreach \j in {0,1,2,3}{\foreach \k in {0,1,2,3}{
  \fill[black!59] (\i,\j,\k) circle[radius=.66pt];
}}}
\foreach \v in {(0,1,0),(0,1,1),(1,1,2),(1,2,2),
                (0,0,1),(0,0,2),(1,0,3),(1,1,3)}
  \fill \v circle[radius=1.03pt];
\coordinate (first-junction) at (0,1,2);
\coordinate (last-junction) at (1,2,3);
\foreach \v in {(0,0,0),(0,1,2),(0,0,3),(1,2,3)}
  \fill \v circle[radius=1.4pt];
\node[path label,left=3pt] at (0,0,0) {$(0,0,0)$};
\node[path label,left=3pt] at (0,0,3) {$(0,0,3)$};
\node[path label] (last-label) at ($(last-junction)+(0mm,4.1mm)$) {$(1,2,3)$};
\node[path label,font=\scriptsize,inner sep=.55pt,left=1.5pt]
  at (first-junction) {$(0,1,2)$};
\draw[black!53,line width=.31pt] (last-junction)--(last-label.south);
\node[path label,right=2pt,text=blue!67!black] at (0,1,.60) {$P$};
\node[path label,above=2pt,text=blue!67!black] at (1,.67,3) {$P$};
\node[path label,right=2pt,text=teal!68!black] at (1,2,2.5) {$S$};
\node[path label,left=3pt,text=red!59!black] at (0,0,1.5) {$\Gamma$};
\end{tikzpicture}
\caption{Two staircases in $I_3^{\Box 3}$ from $(0,0,0)$ to $(1,2,3)$.
Under $\Phi$, the staircase through
$(0,1,2)$ maps to $P\vee S$, and the staircase through $(0,0,3)$ maps
to $\Gamma\vee P$.}
\label{fig:GLMY-staircases}
\end{figure}
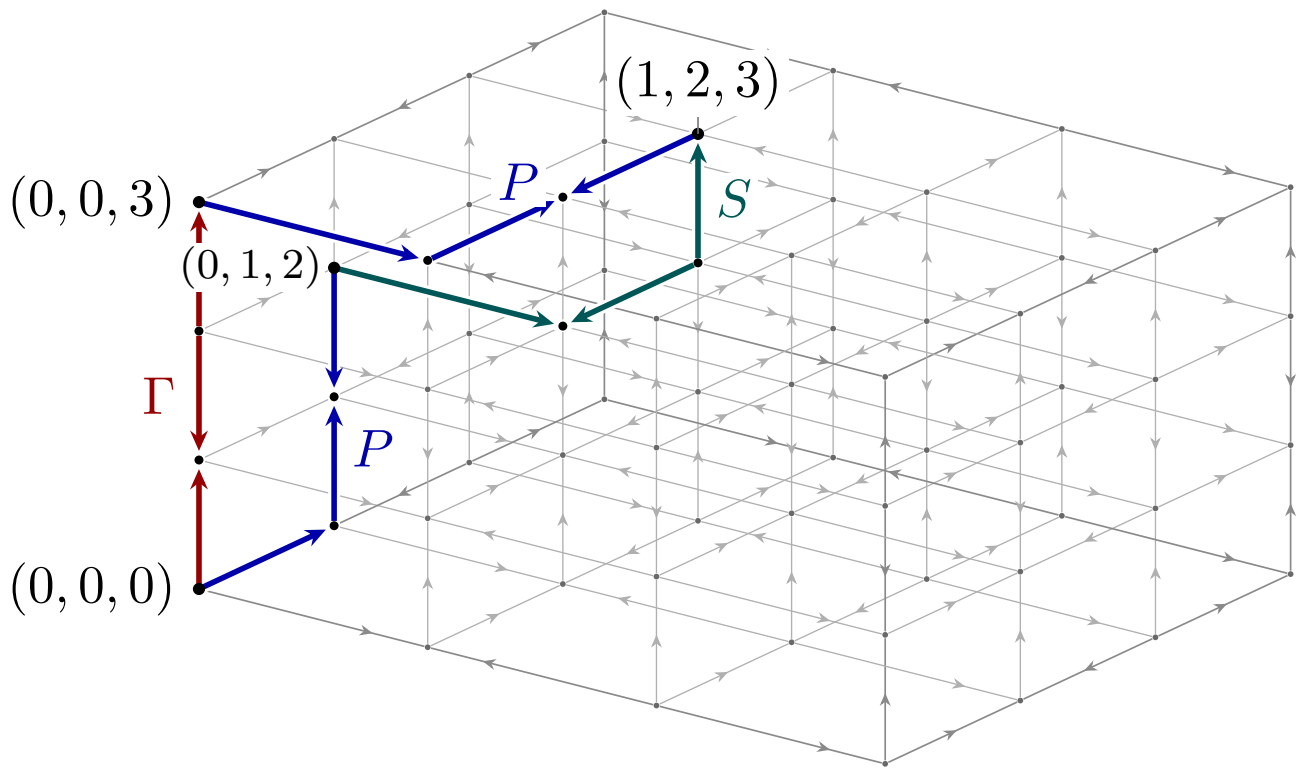

To contract the box with its origin fixed, for $1\le j\le m$ and
$0\le r\le m$ put
\[
 R_{j,r}(t_1,\ldots,t_m)
 =\bigl(\underbrace{0,\ldots,0}_{j-1\ \text{coordinates}},
         \min(t_j,r),t_{j+1},\ldots,t_m\bigr).
\]
Each $R_{j,r}$ is a digraph map. In coordinate $j$ the truncation
preserves the arrows up to vertex $r$ and collapses the arrows
beyond it. The other coordinates are fixed or sent to $0$.
For $1\le r\le m$, the values of $R_{j,r}$ and $R_{j,r-1}$ agree
when $t_j<r$. When $t_j\ge r$, they differ only in coordinate $j$,
where the values are $r$ and $r-1$. Orient the time edge between
these maps according to the arrow between $r$ and $r-1$ in $I_m$.
This gives a pointwise comparison of the same direction at every
vertex. For each $j$ in increasing order, let $r$ decrease from
$m$ to $0$. The last map for coordinate $j$ equals the first map
for coordinate $j+1$, since $R_{j,0}=R_{j+1,m}$ for $j<m$.
The $m^2$ comparisons determine the arrows of the time line
$I_{m^2}$. In this order the maps $R_{j,r}$ are the time slices of
$H:I_m^{\Box m}\Box I_{m^2}\to I_m^{\Box m}$, from
$R_{1,m}=\id$ to the constant map $R_{m,0}$ at the origin.
Every $R_{j,r}$ fixes the origin.

When reversing a path of length $N$, we read its values in reverse order
and reflect its parameter line: $i\to j$ in the reflected line means
$N-i\to N-j$ in the original line. Each of the two full staircases
has length $\binom{m}{2}+m$. Following the first and then the
reverse of the second gives a loop of length $m(m+1)$:
$\lambda:(I_{m(m+1)},\partial I_{m(m+1)})\to(I_m^{\Box m},(0,\ldots,0))$.
Precompose $H$ with $\lambda\Box\id_{I_{m^2}}$ and postcompose with
$\Phi$. The digraph homotopy
\[
 (t,s)\longmapsto \Phi\bigl(H(\lambda(t),s)\bigr)
\]
equals $e_{G,*}$ at $t=0,m(m+1)$ for every time $s$: both endpoints
of $\lambda$ are the origin, and $H$ fixes the origin.
The initial slice is
\[
 \Phi\circ\lambda=P\vee S\vee P^{-}\vee\Gamma^{-}
\]
and the final slice is constant at $e_{G,*}$. On every time edge,
the two slices have the same oriented parameter line and are
pointwise comparable, so identity shrinkings make the comparison
a $C$-homotopy step. The final constant loop is $C$-homotopic to
the length-zero loop by shrinking its parameter line to a point.
Thus $P\vee S\vee P^{-}\vee\Gamma^{-}$ is $C$-null-homotopic.

The contraction of $S$ can be performed inside this concatenation.
Extend each shrinking on the $S$ block by identity maps on the
other blocks using \eqref{eq:GLMY-concatenated-shrinking}. Extend
each pointwise comparison by equality on the other blocks.
The two endpoints of the $S$ block stay at $\rho$, so they still
agree with the endpoints of the adjacent blocks. When $S$ is
replaced by the length-zero loop at $\rho$, the resulting loop is
$P\vee P^{-}\vee\Gamma^{-}$.

To cancel $P\vee P^{-}$, consider the backtrack at its turning point.
Its three successive values have the form $x,y,x$. If $x\ne y$,
the two edges at the turning point use the same arrow of the target,
because the parameter line of the reversed path is reflected.
Replace $y$ by $x$ and leave all other values fixed. This gives a
loop on the same line, since both arrows incident to the changed
vertex collapse and all other arrows are unchanged. The old and
new loops satisfy the pointwise comparison from the new loop to the
old loop when $x\to y$, and from the old loop to the new loop when
$y\to x$. If $x=y$, the loops are equal. Now shrink
the two constant edges to a single vertex with value $x$. This
gives a $C$-homotopic loop with the backtrack removed. Repeating the
cancellation removes $P\vee P^{-}$, and shows that $[\Gamma^{-}]=1$.
The same cancellation gives $[\Gamma][\Gamma^{-}]=1$, so
$[\Gamma]=1$. Since $\Gamma$ was arbitrary, this proves
$\pi_1(L(G,*),e_{G,*})=0$.

All parameter lines and boxes used in the proof are finite, and the
argument applies to infinite target digraphs as well. For $k\ge2$,
apply the conclusion to $L^{k-2}(G,*)$ at its distinguished
basepoint. Its original loop digraph is $L^{k-1}(G,*)$.
By \eqref{eq:GLMY-recursive},
\[
 \pi_k^{\mathrm{GLMY}}(G,*)
 =\pi_1\bigl(L^{k-1}(G,*)\bigr)=0.\qedhere
\]
\end{proof}
\section{Digraph spheres}\label{sec:spheres}

We use boundary quotients of boxes as domains for sphere
representatives. Integral GLMY path homology determines which
quotients are digraph spheres.

\subsection{Integral GLMY path homology}

\begin{definition}\label{def:regular-allowed}
Let $X$ be a digraph and $p\ge0$.  A vertex sequence $e_{x_0\cdots x_p}$ in $X$
is \emph{regular} if $x_{i-1}\ne x_i$ for $1\le i\le p$, and it is
\emph{allowed} if $x_{i-1}\to x_i$ for $1\le i\le p$.
Let $\mathcal R_p(X;\Z)$ be the free abelian group on regular $p$-paths and
$\mathcal A_p(X;\Z)\subseteq\mathcal R_p(X;\Z)$ the subgroup generated by
allowed $p$-paths.  For $p\ge1$, define
\[
 \partial e_{x_0\cdots x_p}
 =\sum_{i=0}^p(-1)^i
 e_{x_0\cdots\widehat{x_i}\cdots x_p},
\]
where every nonregular term is set equal to zero.
This boundary is induced by alternating deletion on all vertex sequences.
For two equal consecutive vertices, the two terms that delete one
of them cancel. Each other deletion leaves two equal consecutive
vertices in the resulting path. The boundary thus preserves the subgroup generated by the
nonregular paths, and the induced boundary satisfies
$\partial^2=0$.
\end{definition}

\begin{definition}\label{def:path-complex}
Let $X$ be a digraph.  For $p\ge1$, set
\[
 \Omega_p(X;\Z)
 =\{u\in\mathcal A_p(X;\Z):\partial u\in\mathcal A_{p-1}(X;\Z)\}.
\]
We use the unaugmented convention
\[
 \Omega_0(X;\Z)=\mathcal A_0(X;\Z),
 \qquad
 \partial_0=0,
\]
and define
\[
 H_p(X;\Z)=H_p(\Omega_*(X;\Z)).
\]
A chain in $\Omega_p(X;\Z)$ is called \emph{$\partial$-invariant}. A cycle is
a $\partial$-invariant chain $u$ satisfying $\partial u=0$.
\end{definition}

A digraph map $f:X\to G$ induces a chain map
$f_\#:\Omega_*(X;\Z)\to\Omega_*(G;\Z)$ by applying $f$ to each
vertex of a regular path and setting every image with consecutive repeated
vertices equal to zero.  We write
$H_p(f):H_p(X;\Z)\to H_p(G;\Z)$ for the induced map in degree $p$.
Homotopic digraph maps induce the same map in GLMY path homology
\cite[Theorem~3.3]{GLMY}. For infinite digraphs, restrict the source to the
finite induced subdigraph on the vertices of a cycle. The homotopy on
this subdigraph has image in a finite induced subdigraph of the target.

\subsection{Sphere presentations and standard models}

\begin{definition}\label{def:box-presentation}
Let $k\ge1$ and let
\[
 \mathcal I=(I_{m_1},\ldots,I_{m_k})
\]
be a $k$-tuple of finite oriented lines of positive lengths. The choice
of $I_{m_i}$ includes its orientation, which may differ between factors
of the same length.  Its \emph{box}
and \emph{full coordinate boundary} are
\[
 B_{\mathcal I}=I_{m_1}\Box\cdots\Box I_{m_k},
\]
and
\[
 V_{\partial B_{\mathcal I}}
 =\{x\in V_{B_{\mathcal I}}:x_i\in\{0,m_i\}
   \text{ for some }i\},
\]
where the boundary is the induced subdigraph on the displayed vertex set.
Let $q_{\mathcal I}:V_{B_{\mathcal I}}\to
(V_{B_{\mathcal I}}\setminus V_{\partial B_{\mathcal I}})\sqcup\{*\}$ be the
vertex map that is the identity off the boundary and sends every boundary
vertex to $*$.  The quotient digraph
$S_{\mathcal I}=B_{\mathcal I}/\partial B_{\mathcal I}$ has vertex set
\[
 V_{S_{\mathcal I}}
 =(V_{B_{\mathcal I}}\setminus V_{\partial B_{\mathcal I}})\sqcup\{*\}.
\]
For distinct vertices $x,y\in V_{S_{\mathcal I}}$, there is an arrow
$x\to y$ if and only if there are vertices $u,v\in V_{B_{\mathcal I}}$ with
$u\to v$, $q_{\mathcal I}(u)=x$, and $q_{\mathcal I}(v)=y$.  Thus arrows
whose endpoints are both collapsed to $*$ are discarded. The quotient
map gives a digraph map of pairs
\[
 q_{\mathcal I}:(B_{\mathcal I},\partial B_{\mathcal I})
 \longrightarrow(S_{\mathcal I},*).
\]
\end{definition}

\begin{definition}\label{def:digraph-sphere}
Let $\mathcal I$ be a $k$-tuple as in \cref{def:box-presentation}.
The box quotient $(S_{\mathcal I},*)$ is a \emph{$k$-dimensional digraph
sphere} if
\[
 H_q(S_{\mathcal I};\Z)\cong
 \begin{cases}
  \Z,&q=0,k,\\
  0,&q\notin\{0,k\}.
 \end{cases}
\]
Let $\mathfrak S_k$ be the collection of all $k$-tuples $\mathcal I$ whose
boundary quotients $(S_{\mathcal I},*)$ are $k$-dimensional digraph spheres.
For $\mathcal I\in\mathfrak S_k$, the quotient
digraph map $q_{\mathcal I}$ is called a \emph{$k$-sphere presentation}.
\end{definition}

Every $\mathcal I=(I_{m_1},\ldots,I_{m_k})\in\mathfrak S_k$ satisfies
$m_i\ge3$ for all $i$. If $m_i=1$, the whole box is its boundary.
If $m_i=2$, every interior vertex has $i$th coordinate~$1$.
Replacing this coordinate by~$0$ gives a boundary vertex. The orientation
of the edge between $0$ and $1$ in $I_{m_i}$ therefore gives either
$*\to x$ for every interior vertex $x$, or $x\to *$ for every such $x$.
The identity of $S_{\mathcal I}$ is then homotopic to the constant map at
$*$ in one step, so its positive-degree GLMY path homology vanishes.

\begin{definition}\label{def:sphere-representative}
Let $(G,*)$ be a based digraph and $q_{\mathcal I}$ a $k$-sphere presentation.
A \emph{$k$-sphere
representative in $(G,*)$ on $q_{\mathcal I}$} is a digraph map of pairs
\[
 f:(B_{\mathcal I},\partial B_{\mathcal I})\longrightarrow(G,*).
\]
Equivalently, $f$ factors uniquely through a based digraph map
\[
 \bar f:(S_{\mathcal I},*)\longrightarrow(G,*),
 \qquad
 \bar f(q_{\mathcal I}(x))=f(x)\quad(x\in V_{B_{\mathcal I}}).
\]
The formula is well defined because $f$ has value $*$ on the boundary.
Each arrow of $S_{\mathcal I}$ lifts to an arrow of $B_{\mathcal I}$,
so $\bar f$ is a digraph map.
\end{definition}

\begin{definition}\label{def:standard-line}
For $m\ge0$, let $J_m$ be the alternating oriented line on
$\{0,\ldots,m\}$ in which each even vertex points to its adjacent odd
vertices. For $n\ge1$, the \emph{standard alternating line} of length
$2n$ is
\[
 J_{2n}=[0\to1\leftarrow2\to\cdots\to2n-1\leftarrow2n].
\]
The standard line of length zero is the one-vertex digraph $J_0$.
\end{definition}

Let $k\ge1$ and $\nn=(n_1,\ldots,n_k)$ with $n_i\ge2$.
Applying \cref{def:box-presentation} to the tuple
$(J_{2n_1},\ldots,J_{2n_k})$, we write
\[
 \begin{aligned}
 &B_{\nn}=J_{2n_1}\Box\cdots\Box J_{2n_k},\\
 &V_{\partial B_{\nn}}=\{x\in V_{B_{\nn}}:x_i\in\{0,2n_i\}
 \text{ for some }i\},\\
 &\Sn=B_{\nn}/\partial B_{\nn}.
 \end{aligned}
\]
Write $q_{\nn}$ for the quotient map of \cref{def:box-presentation}, so
\[
 q_{\nn}:(B_{\nn},\partial B_{\nn})\longrightarrow(\Sn,*).
\]
We call $B_{\nn}$ a \emph{standard alternating box}.

\begin{theorem}\label{thm:admissibility}
For every $\nn=(n_1,\ldots,n_k)$ with $n_i\ge2$, the quotient
$(\Sn,*)$ is a $k$-dimensional digraph sphere.
\end{theorem}

Appendix~\ref{app:homology} computes $H_q(\Sn;\Z)$ using the product
of the interior lines relative to their coordinate boundary.

\begin{definition}\label{def:standard-representative}
For $\nn$ as in \cref{thm:admissibility}, we call $(\Sn,*)$ a
\emph{standard digraph sphere} and the quotient digraph map $q_{\nn}$ a
\emph{standard $k$-sphere presentation}.
A sphere representative on $q_{\nn}$ is called a \emph{standard representative}.
\end{definition}

For example,
\[
 S^2_{4\times4}=\bigl(J_4\Box J_4\bigr)/\partial\bigl(J_4\Box J_4\bigr)
\]
has the nine vertices $(i,j)$, $1\le i,j\le3$, together with the collapsed
boundary vertex $*$.

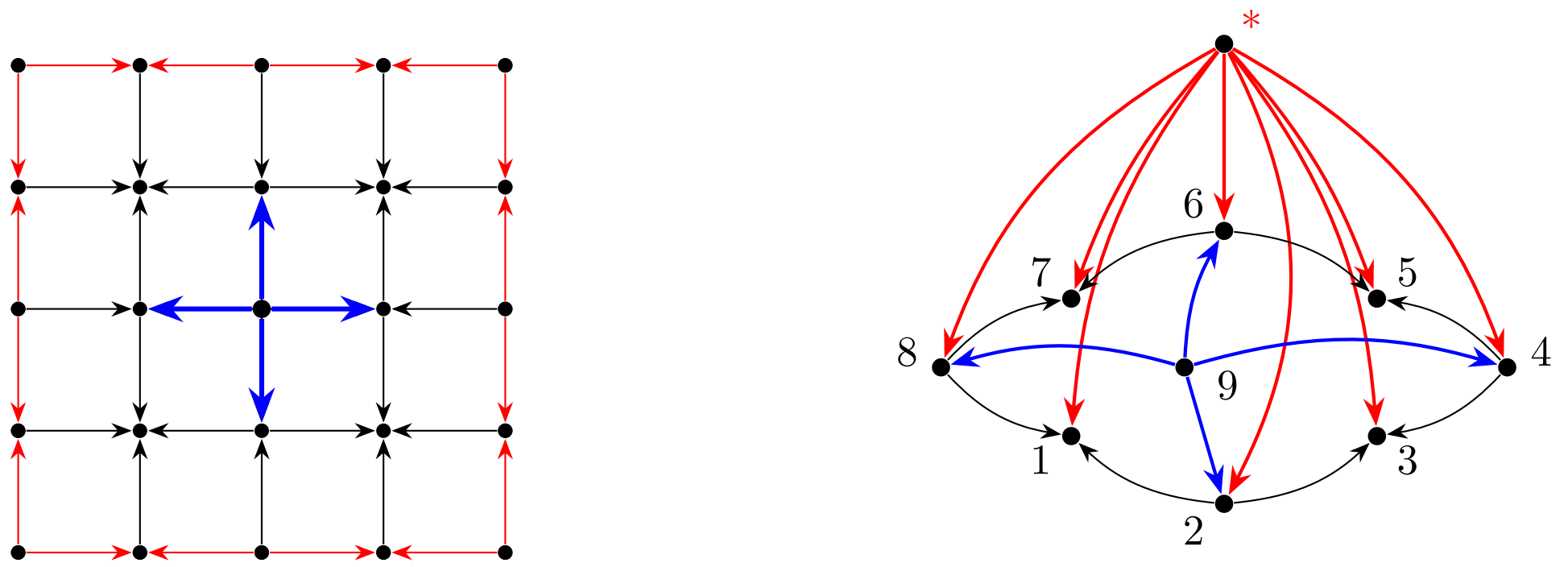
\begin{figure}[H]
\centering
\begin{minipage}[t]{0.48\textwidth}
\centering
\begin{tikzpicture}[scale=1,>=Stealth]
\foreach \i in {0,...,4}{\foreach \j in {0,...,4}{\node[small graph vertex] (p\i\j) at (\i,\j) {};}}
\foreach \i in {0,...,4}{\node[small graph vertex] at (\i,0) {};\node[small graph vertex] at (\i,4) {};}
\foreach \j in {1,2,3}{\node[small graph vertex] at (0,\j) {};\node[small graph vertex] at (4,\j) {};}
\foreach \j in {1,2,3}{\draw[->,thin] (p0\j)--(p1\j);\draw[->,thin] (p2\j)--(p1\j);\draw[->,thin] (p2\j)--(p3\j);\draw[->,thin] (p4\j)--(p3\j);}
\foreach \j in {0,4}{\draw[->,red,thin] (p0\j)--(p1\j);\draw[->,red,thin] (p2\j)--(p1\j);\draw[->,red,thin] (p2\j)--(p3\j);\draw[->,red,thin] (p4\j)--(p3\j);}
\foreach \i in {1,2,3}{\draw[->,thin] (p\i0)--(p\i1);\draw[->,thin] (p\i2)--(p\i1);\draw[->,thin] (p\i2)--(p\i3);\draw[->,thin] (p\i4)--(p\i3);}
\foreach \i in {0,4}{\draw[->,red,thin] (p\i0)--(p\i1);\draw[->,red,thin] (p\i2)--(p\i1);\draw[->,red,thin] (p\i2)--(p\i3);\draw[->,red,thin] (p\i4)--(p\i3);}
\node[graph vertex] (p22) at (2,2) {};
\draw[highlighted edge,blue] (p22)--(p12);
\draw[highlighted edge,blue] (p22)--(p32);
\draw[highlighted edge,blue] (p22)--(p21);
\draw[highlighted edge,blue] (p22)--(p23);
\end{tikzpicture}

\smallskip
\small $J_4\Box J_4$, with the full coordinate boundary in red.
\end{minipage}
\hfill
\begin{minipage}[t]{0.48\textwidth}
\centering
\begin{tikzpicture}[xscale=.93,yscale=.83,>=Stealth]
\node[graph vertex] (s) at (0,5.2) {};
\node[font=\small,text=red] at (.24,5.44) {$*$};
\node[graph vertex] (x8) at (-2.5,2.0) {};
\node[graph vertex] (x4) at ( 2.5,2.0) {};
\node[graph vertex] (x6) at (0,3.35) {};
\node[graph vertex] (x2) at (0,.65) {};
\node[graph vertex] (b1) at (-1.35,2.68) {};
\node[graph vertex] (a1) at ( 1.35,2.68) {};
\node[graph vertex] (b2) at (-1.35,1.32) {};
\node[graph vertex] (a2) at ( 1.35,1.32) {};
\node[graph vertex] (x9) at (-0.35,2.0) {};
\draw[->,red,thick,bend right=20] (s) to (x8);
\draw[->,red,thick,bend left=20] (s) to (x4);
\draw[->,red,thick] (s) -- (x6);
\draw[->,red,thick,bend left=25] (s) to (x2);
\draw[->,red,thick,bend right=10] (s) to (b1);
\draw[->,red,thick,bend left=10] (s) to (a1);
\draw[->,red,thick,bend right=15] (s) to (b2);
\draw[->,red,thick,bend left=15] (s) to (a2);
\draw[->,blue,thick,bend left=18] (x9) to (x4);
\draw[->,blue,thick,bend right=18] (x9) to (x8);
\draw[->,blue,thick] (x9) -- (x2);
\draw[->,blue,thick,bend left=12] (x9) to (x6);
\draw[->,bend left=20] (x8) to (b1);
\draw[->,bend right=20] (x6) to (b1);
\draw[->,bend left=20] (x6) to (a1);
\draw[->,bend right=20] (x4) to (a1);
\draw[->,bend right=20] (x8) to (b2);
\draw[->,bend left=20] (x2) to (b2);
\draw[->,bend right=20] (x2) to (a2);
\draw[->,bend left=20] (x4) to (a2);
\node[font=\small] at (-2.80,2.15) {$8$};
\node[font=\small] at (-1.62,2.94) {$7$};
\node[font=\small] at (-.27,3.62) {$6$};
\node[font=\small] at ( 1.62,2.94) {$5$};
\node[font=\small] at ( 2.80,2.15) {$4$};
\node[font=\small] at ( 1.62,1.08) {$3$};
\node[font=\small] at (-1.62,1.08) {$1$};
\node[font=\small] at (-.27,.38) {$2$};
\node[font=\small] at (.03,1.82) {$9$};
\end{tikzpicture}

\smallskip
\small $(S^2_{4\times4},*)$ after the boundary collapse.
\end{minipage}
\caption{The standard digraph sphere $(S^2_{4\times4},*)$.  Left: the standard alternating box $J_4\Box J_4$, with the full coordinate boundary indicated in red and the arrows from the central vertex $(2,2)$ highlighted in blue.  Right: the quotient obtained by collapsing the full coordinate boundary to the basepoint~$*$.
The labels $1,\ldots,8$ correspond, in order, to
$(1,1),(2,1),(3,1),(3,2),(3,3),(2,3),(1,3),(1,2)$; label $9$ denotes $(2,2)$.}
\label{fig:standard-sphere-4x4}
\end{figure}

\section{\texorpdfstring{$C$}{C}-homotopy and standardization}\label{sec:shrinking}

We compare sphere representatives on different boxes by pulling
them back to a common box. For standard boxes, we show that terminal
extensions and homotopies relative to the boundary suffice.

\subsection{Shrinkings and sphere representatives}

\begin{definition}\label{def:line-extension}
For a based loop $\gamma$ on $I_m$ and a shrinking $h:I_N\to I_m$,
the \emph{extension of $\gamma$ along $h$} is the based loop
$\gamma\circ h$ on $I_N$.
\end{definition}

\begin{definition}\label{def:shrinking}
Let $\mathcal I=(I_{m_1},\ldots,I_{m_k})$ and
$\mathcal J=(I_{N_1},\ldots,I_{N_k})$ be $k$-tuples of finite oriented
lines of positive lengths, and write
\[
 B_{\mathcal J}=I_{N_1}\Box\cdots\Box I_{N_k},
 \qquad
 B_{\mathcal I}=I_{m_1}\Box\cdots\Box I_{m_k}.
\]
A \emph{product shrinking} is a digraph map
\[
 \rho=\rho_1\Box\cdots\Box\rho_k:
 B_{\mathcal J}\longrightarrow B_{\mathcal I}
\]
for which every coordinate digraph map
$\rho_i:I_{N_i}\to I_{m_i}$ is a shrinking.  A product shrinking
sends the full coordinate boundary into the full coordinate boundary.
\end{definition}

\begin{definition}\label{def:product-extension}
Let
\[
 f:(B_{\mathcal I},\partial B_{\mathcal I})\longrightarrow(G,*)
\]
be a sphere representative, let $\mathcal J\in\mathfrak S_k$, and let
$\rho:B_{\mathcal J}\to B_{\mathcal I}$ be a product shrinking.  The sphere
representative
\[
 f\circ\rho:(B_{\mathcal J},\partial B_{\mathcal J})\longrightarrow(G,*)
\]
is called the \emph{extension of $f$ along $\rho$}.
\end{definition}

\begin{definition}\label{def:direct-sphere-C}
Let
\[
 f:(B_{\mathcal I},\partial B_{\mathcal I})\to(G,*),
 \qquad
 g:(B_{\mathcal J},\partial B_{\mathcal J})\to(G,*)
\]
be $k$-sphere representatives.  A \emph{one-step direct $C$-homotopy from
$f$ to $g$} consists of $\mathcal K\in\mathfrak S_k$ and product shrinkings
\[
 \rho:B_{\mathcal K}\to B_{\mathcal I},
 \qquad
 \sigma:B_{\mathcal K}\to B_{\mathcal J}
\]
such that
\[
 f\rho(u)\toeq g\sigma(u)
 \qquad\text{for every }u\in V_{B_{\mathcal K}}.
\]
A \emph{one-step inverse $C$-homotopy from $f$ to $g$} is a one-step direct
$C$-homotopy from $g$ to $f$.
\end{definition}

\begin{definition}\label{def:sphere-C}
Let $(G,*)$ be a based digraph and $k\ge1$.
Two $k$-sphere representatives $f$ and $g$ in $(G,*)$ are \emph{$C$-homotopic}, written
$f\simeq_C g$, if there is a finite sequence
\[
 f=f_0,f_1,\ldots,f_r=g
\]
such that, for $0\le j<r$, $f_j$ and $f_{j+1}$ are related by a one-step direct
or one-step inverse $C$-homotopy.  The set of $k$-dimensional $C$-homotopy
classes in $(G,*)$ is
\begin{equation}\label{eq:piC-definition}
 \pC_k(G,*)=
 \left(
 \bigsqcup_{\mathcal I\in\mathfrak S_k}
 \Hom((B_{\mathcal I},\partial B_{\mathcal I}),(G,*))
 \right)\Bigg/\simeq_C.
\end{equation}
\end{definition}

On a fixed presentation, the identity shrinkings turn each one-step
relative homotopy into a one-step direct or inverse $C$-homotopy.
Thus relative homotopy implies $C$-homotopy.

For an extension along $\rho:B_{\mathcal J}\to B_{\mathcal I}$,
the shrinkings $\rho$ and $\id_{B_{\mathcal J}}$ give a one-step
direct $C$-homotopy between $f$ and $f\circ\rho$, since
\begin{equation}\label{eq:extension-C}
 f(\rho(u))=(f\circ\rho)(u)
 \qquad(u\in V_{B_{\mathcal J}}).
\end{equation}

\subsection{Reduction to standard representatives}

\begin{lemma}\label{lem:standard-cover}
Let $I_m$ be a finite oriented line.  For every integer $N\ge m$ there is a shrinking
\[
 r:J_{2N}\longrightarrow I_m.
\]
\end{lemma}
\begin{proof}
If $m=0$, take the constant map.  For $m\ge1$ and $N\ge m$, define
\[
 \begin{aligned}
 r(2i)&=i &&(0\le i\le m),\\
 r(2i+1)&=
 \begin{cases}
 i+1,&i\to i+1\text{ in }I_m,\\
 i,&i+1\to i\text{ in }I_m,
 \end{cases} &&(0\le i<m),\\
 r(t)&=m &&(2m\le t\le2N).
 \end{aligned}
\]
On each pair of arrows $2i\to2i+1\leftarrow2i+2$, the values are
\[
 (r(2i),r(2i+1),r(2i+2))=
 \begin{cases}
 (i,i+1,i+1),&i\to i+1\text{ in }I_m,\\
 (i,i,i+1),&i+1\to i\text{ in }I_m.
 \end{cases}
\]
In each pair, one arrow maps to the target arrow and the other
collapses. All arrows after $2m$ collapse. The values increase
weakly from $0$ to $m$ and include every integer in this range, so
$r$ is a shrinking.
\end{proof}

For example, the target line $I_3=[0\leftarrow1\to2\to3]$ is
obtained from $J_6$ by the shrinking with vertex-value list
\[
 r=(0,0,1,2,2,3,3):J_6\longrightarrow I_3.
\]
Each repeated value marks a collapsed edge.

\begin{lemma}\label{lem:common-refinement}
Let $I_{m_1},\ldots,I_{m_r}$ and $I_m$ be finite oriented lines, with $r\ge1$, and let
\[
 h_\nu:I_{m_\nu}\longrightarrow I_m
 \qquad(1\le \nu\le r)
\]
be shrinkings.  Then there are a standard alternating line $J_{2N}$ and shrinkings
\[
 u_\nu:J_{2N}\longrightarrow I_{m_\nu}
 \qquad(1\le \nu\le r)
\]
such that
\begin{equation}\label{eq:common-refinement-equality}
 h_1u_1=h_2u_2=\cdots=h_ru_r:J_{2N}\longrightarrow I_m.
\end{equation}
\end{lemma}

\begin{proof}
For $1\le\nu\le r$ and $0\le j\le m$, define the endpoints of each fiber by
\[
 a_{\nu,j}=\min h_\nu^{-1}(j),\qquad
 b_{\nu,j}=\max h_\nu^{-1}(j).
\]
Monotonicity and surjectivity give
\[
 h_\nu^{-1}(j)=\{a_{\nu,j},\ldots,b_{\nu,j}\},
\]
and
\begin{equation}\label{eq:fiber-endpoints}
 a_{\nu,0}=0,\qquad b_{\nu,m}=m_\nu,\qquad
 a_{\nu,j+1}=b_{\nu,j}+1\quad(0\le j<m).
\end{equation}
In $h_1^{-1}(j)\times\cdots\times h_r^{-1}(j)$, start at
$(a_{1,j},\ldots,a_{r,j})$ and increase each coordinate in turn to
its upper endpoint. When coordinate $\nu$ increases from $s$ to $s+1$,
orient the parameter edge as in $I_{m_\nu}$. Its projection to that
factor preserves the arrow, and every other projection collapses it.

Between the list over $j$ and the list over $j+1$, insert one parameter
edge joining
\[
 (b_{1,j},\ldots,b_{r,j})
 \quad\text{and}\quad
 (a_{1,j+1},\ldots,a_{r,j+1}).
\]
Since $h_\nu(b_{\nu,j})=j$ and $h_\nu(a_{\nu,j+1})=j+1$,
\begin{equation}\label{eq:transition-orientation}
 \begin{cases}
 b_{\nu,j}\to a_{\nu,j+1},&j\to j+1\text{ in }I_m,\\
 a_{\nu,j+1}\to b_{\nu,j},&j+1\to j\text{ in }I_m.
 \end{cases}
\end{equation}
Give this edge of the parameter line the direction of the edge
between $j$ and $j+1$ in $I_m$. Every projection then preserves its
orientation, even though all coordinates change together.

The concatenated list has
\[
 M=m+\sum_{j=0}^m\sum_{\nu=1}^r(b_{\nu,j}-a_{\nu,j})
   =\sum_{\nu=1}^r m_\nu-(r-1)m
\]
edges.  Identify its ordered vertices with $0,\ldots,M$ and call the
resulting oriented line $I_M$.  Its coordinate projections
$\pi_\nu:I_M\to I_{m_\nu}$ satisfy
\[
 \begin{gathered}
 \pi_\nu(0)=0,\qquad \pi_\nu(M)=m_\nu,\qquad
 \pi_\nu(t+1)-\pi_\nu(t)\in\{0,1\}\quad(0\le t<M),\\
 \pi_\nu(V_{I_M})=\bigcup_{j=0}^m h_\nu^{-1}(j)=V_{I_{m_\nu}}.
 \end{gathered}
\]
The chosen orientations make each projection a digraph map, so
each $\pi_\nu$ is a shrinking. At a vertex $t$ over $j$,
$h_\nu(\pi_\nu(t))=j$ for every $\nu$. Hence
\begin{equation}\label{eq:IM-common-composite}
 h_1\pi_1=\cdots=h_r\pi_r:I_M\longrightarrow I_m.
\end{equation}
Choose $N\ge M$ and a shrinking $s:J_{2N}\to I_M$ by
\cref{lem:standard-cover}.  The shrinkings $u_\nu=\pi_\nu s$ satisfy
\[
 h_1u_1=(h_1\pi_1)s=\cdots=(h_r\pi_r)s=h_ru_r,
\]
as required.  The construction includes $m=0$ and singleton fibers.
\end{proof}

\begin{figure}[H]
\centering
\begin{tikzpicture}[
  x=1cm,y=1cm,>=Stealth,
  every node/.style={font=\small},
  vertex/.style={circle,fill=black,inner sep=1.5pt,outer sep=0pt},
  edge/.style={->,line width=.45pt,shorten <=1pt,shorten >=1pt},
  background edge/.style={edge,black!45},
  selected edge/.style={edge,black,line width=1pt},
  map arrow/.style={->,line width=.6pt}
]
\node[anchor=west] at (0,1.25) {(a) The two shrinkings};
\node at (2.0,.70) {$I_3$};
\node at (7.0,.70) {$I_2$};
\node at (12.5,.70) {$I_4$};
\foreach \x/\i/\value in {.5/0/0,1.5/1/1,2.5/2/1,3.5/3/2}{
  \node[vertex] (a\i) at (\x,0) {};
  \node[above=4pt] at (a\i) {$\i$};
  \node[below=6pt] at (a\i) {$\value$};
}
\draw[edge] (a0)--(a1);
\draw[edge] (a1)--(a2);
\draw[edge] (a3)--(a2);
\foreach \x/\i in {6/0,7/1,8/2}{
  \node[vertex] (t\i) at (\x,0) {};
  \node[above=4pt] at (t\i) {$\i$};
}
\draw[edge] (t0)--(t1);
\draw[edge] (t2)--(t1);
\foreach \x/\i/\value in {10.5/0/0,11.5/1/1,12.5/2/1,13.5/3/1,14.5/4/2}{
  \node[vertex] (b\i) at (\x,0) {};
  \node[above=4pt] at (b\i) {$\i$};
  \node[below=6pt] at (b\i) {$\value$};
}
\draw[edge] (b0)--(b1);
\draw[edge] (b2)--(b1);
\draw[edge] (b2)--(b3);
\draw[edge] (b4)--(b3);
\draw[map arrow] (4.15,0)--node[above=4pt] {$h_1$} (5.35,0);
\draw[map arrow] (9.85,0)--node[above=4pt] {$h_2$} (8.65,0);

\node[anchor=west] at (0,-1.35)
  {(b) The fiber $h_1^{-1}(1)\times h_2^{-1}(1)$};
\foreach \x/\i in {1.55/1,3.25/2}{
  \foreach \y/\j in {-4.15/1,-3.25/2,-2.35/3}{
    \node[vertex] (r\i\j) at (\x,\y) {};
  }
}
\foreach \j in {1,2,3}{
  \node[left=5pt] at (r1\j) {$(1,\j)$};
  \node[right=5pt] at (r2\j) {$(2,\j)$};
}
\draw[background edge] (r12)--(r11);
\draw[background edge] (r12)--(r13);
\draw[background edge] (r12)--(r22);
\draw[background edge] (r13)--(r23);
\draw[selected edge] (r11)--(r21);
\draw[selected edge] (r22)--(r21);
\draw[selected edge] (r22)--(r23);

\node[anchor=west] at (6.1,-1.35) {(c) The common refinement};
\node[anchor=east] at (7.55,-2.30) {$I_5$};
\node[anchor=east] at (7.55,-3.00) {$\pi_1$};
\node[anchor=east] at (7.55,-3.60) {$\pi_2$};
\foreach \x/\i/\first/\second in {
  8.3/0/0/0,9.5/1/1/1,10.7/2/2/1,
  11.9/3/2/2,13.1/4/2/3,14.3/5/3/4}{
  \node[vertex] (c\i) at (\x,-2.30) {};
  \node[above=4pt] at (c\i) {$\i$};
  \node at (\x,-3.00) {$\first$};
  \node at (\x,-3.60) {$\second$};
}
\draw[edge] (c0)--(c1);
\draw[edge] (c1)--(c2);
\draw[edge] (c3)--(c2);
\draw[edge] (c3)--(c4);
\draw[edge] (c5)--(c4);
\node at (10.50,-4.40) {$h_1\pi_1=h_2\pi_2=(0,1,1,1,1,2)$};
\end{tikzpicture}
\caption{A common refinement of $h_1:I_3\to I_2$ and $h_2:I_4\to I_2$.
(a) Lower labels give the images under the shrinkings.
(b) The thick staircase is read in increasing coordinate order; its arrows
retain their coordinate directions.
(c) Joining the staircase to the singleton fibers over $0$ and $2$ gives
$I_5$, with the two shrinking projections displayed vertex by vertex.}
\label{fig:common-refinement-r2}
\end{figure}
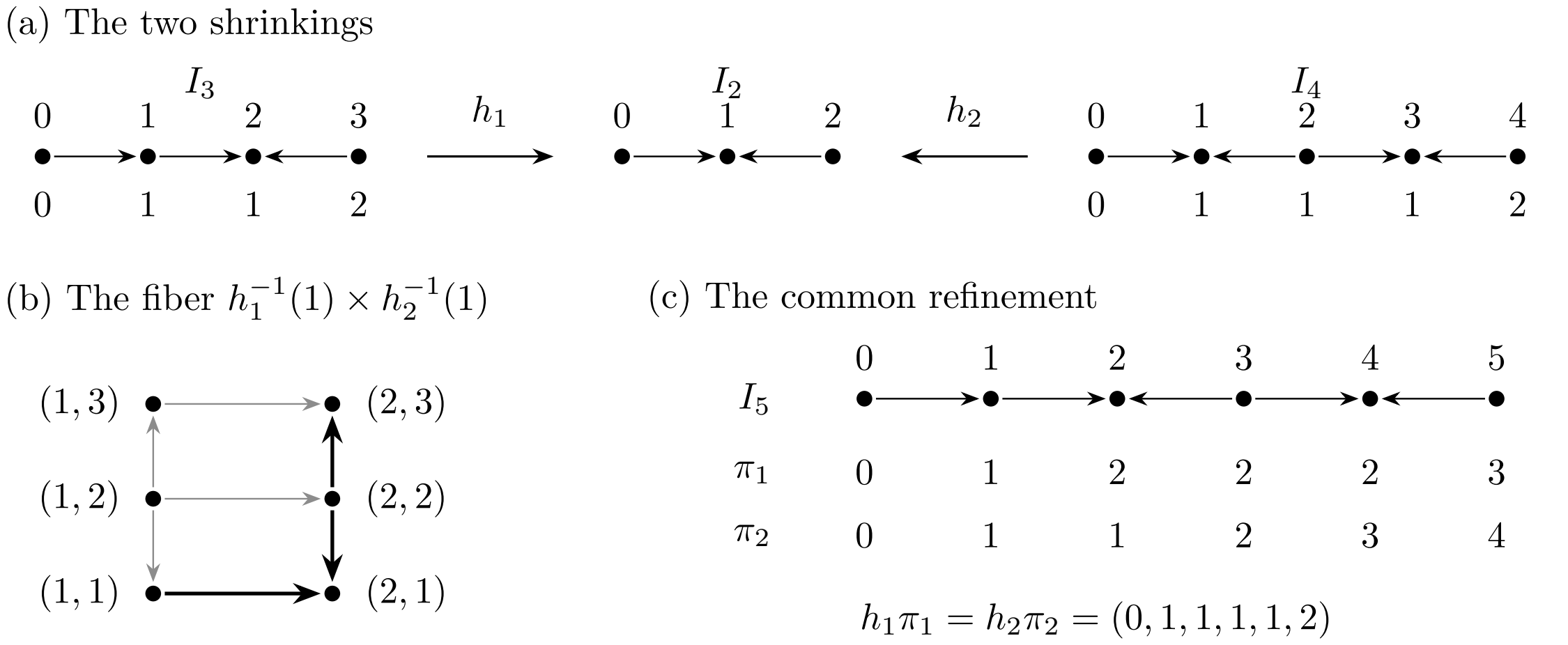

Figure~\ref{fig:common-refinement-r2} gives an example. The fiber over
$1$ is traversed one coordinate at a time. At either end of this fiber,
the two coordinates cross the next target edge together. The resulting
projections have the same composite with $h_1$ and $h_2$.

\Needspace{7\baselineskip}
\begin{proposition}[Reduction to standard representatives]\label{prop:standard-reduction}
Let $(G,*)$ be a based digraph and $k\ge1$.
\begin{enumerate}
\item Every sphere representative
\[
 f:(B_{\mathcal I},\partial B_{\mathcal I})\longrightarrow(G,*),
 \qquad \mathcal I\in\mathfrak S_k,
\]
is $C$-homotopic to a standard representative.
\item If two standard representatives are $C$-homotopic through arbitrary
sphere presentations, then the $C$-homotopy can be chosen to use only
presentations of standard digraph spheres.
\end{enumerate}
Consequently the inclusion of the standard representatives induces a natural
bijection
\begin{equation}\label{eq:standard-reduction}
 \pC_k(G,*)\cong
 \left(
 \bigsqcup_{n_1,\ldots,n_k\ge2}
 \Hom\!\left((B_{\nn},\partial B_{\nn}),(G,*)\right)
 \right)\Bigg/\simeq_C^{\mathrm{std}},
\end{equation}
where $B_{\nn}=J_{2n_1}\Box\cdots\Box J_{2n_k}$
and $\simeq_C^{\mathrm{std}}$ is generated by one-step $C$-homotopies
whose representatives and common sources are all standard.
\end{proposition}

\begin{proof}
For (i), let $\mathcal I=(I_{m_1},\ldots,I_{m_k})\in\mathfrak S_k$.
Choose $N_i\ge\max\{2,m_i\}$ and shrinkings
$\rho_i:J_{2N_i}\to I_{m_i}$ by \cref{lem:standard-cover}, and set
\[
 \rho=\rho_1\Box\cdots\Box\rho_k:
 B_{\mathbf N}\longrightarrow B_{\mathcal I},
 \qquad \mathbf N=(N_1,\ldots,N_k).
\]
Since every $N_i\ge2$, the quotient
$B_{\mathbf N}/\partial B_{\mathbf N}$ is a standard digraph sphere by
\cref{thm:admissibility}. Thus $f\rho$ is a standard representative
and $f\simeq_C f\rho$ by \eqref{eq:extension-C}.

For (ii), write the given $C$-homotopy as
$f=f_0,f_1,\ldots,f_r=g$, where
$f_j:(B_{\mathcal I_j},\partial B_{\mathcal I_j})\to(G,*)$
and consecutive terms are related by a one-step direct or inverse
$C$-homotopy.  There is nothing to prove if $r=0$.
For $0\le j<r$, choose $\mathcal K_j\in\mathfrak S_k$ and product
shrinkings witnessing the $j$th step:
\[
 a_j:B_{\mathcal K_j}\to B_{\mathcal I_j},
 \qquad
 b_j:B_{\mathcal K_j}\to B_{\mathcal I_{j+1}}.
\]
Thus
\begin{equation}\label{eq:std-step}
 f_ja_j\toeq f_{j+1}b_j
 \quad\text{or}\quad
 f_{j+1}b_j\toeq f_ja_j.
\end{equation}
As in (i), choose a product shrinking
$\tau_j:B_{\mathbf N_j}\to B_{\mathcal K_j}$ with every
$N_{j,i}\ge2$. Precomposing \eqref{eq:std-step} with $\tau_j$ gives
\[
 f_ja_j\tau_j\toeq f_{j+1}b_j\tau_j
 \quad\text{or}\quad
 f_{j+1}b_j\tau_j\toeq f_ja_j\tau_j
\]
on the standard box $B_{\mathbf N_j}$.

At $f_j$, for $1\le j<r$, the previous and next steps use the
extensions $f_jb_{j-1}\tau_{j-1}$ and $f_ja_j\tau_j$, respectively.
Both product shrinkings $b_{j-1}\tau_{j-1}$ and $a_j\tau_j$ have
target $B_{\mathcal I_j}$. Applying \cref{lem:common-refinement}
to their coordinate maps gives product
shrinkings
\[
 \rho_j:B_{\mathbf M_j}\to B_{\mathbf N_{j-1}},
 \qquad
 \sigma_j:B_{\mathbf M_j}\to B_{\mathbf N_j},
\]
with $\rho_j=\rho_{j,1}\Box\cdots\Box\rho_{j,k}$ and
$\sigma_j=\sigma_{j,1}\Box\cdots\Box\sigma_{j,k}$, such that
\begin{equation}\label{eq:coordinate-refinement}
 (b_{j-1}\tau_{j-1})_i\rho_{j,i}
   =(a_j\tau_j)_i\sigma_{j,i}
 \qquad(1\le i\le k).
\end{equation}
Taking the product over $i$ gives
\begin{equation}\label{eq:product-refinement}
 b_{j-1}\tau_{j-1}\rho_j=a_j\tau_j\sigma_j.
\end{equation}
Surjectivity implies
$M_{j,i}\ge\max\{N_{j-1,i},N_{j,i}\}\ge2$.
Thus $B_{\mathbf M_j}/\partial B_{\mathbf M_j}$ is a standard digraph
sphere by \cref{thm:admissibility}.  After composition with $f_j$, \eqref{eq:product-refinement} gives
\[
 (f_jb_{j-1}\tau_{j-1})\rho_j=(f_ja_j\tau_j)\sigma_j.
\]
This equality gives a one-step direct $C$-homotopy on the common
standard box.

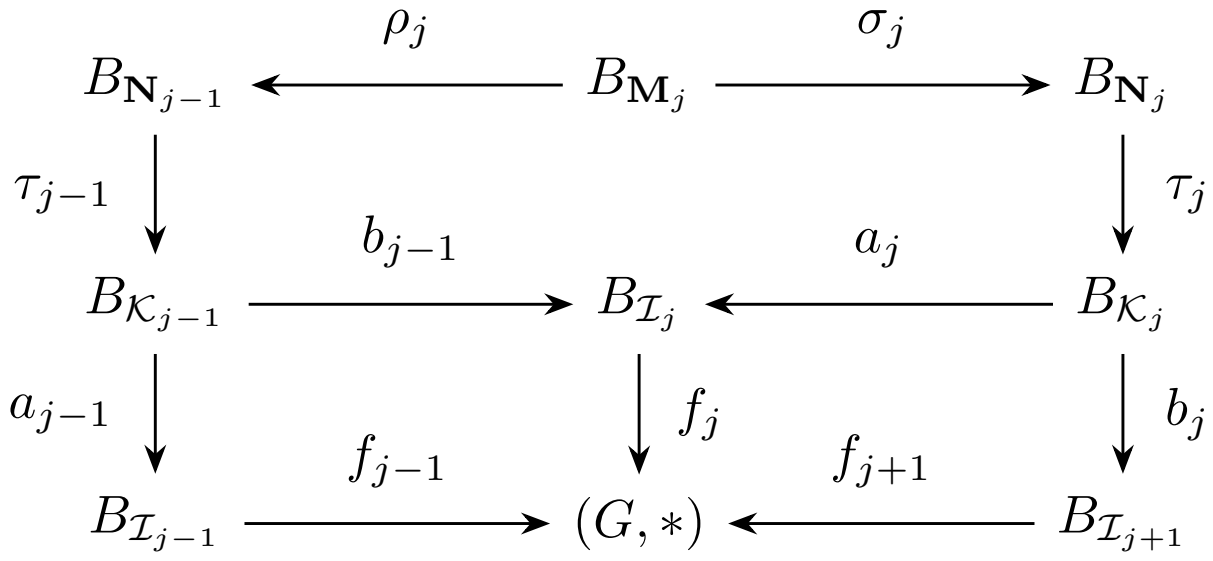
\begin{figure}[H]
\centering
\begin{tikzpicture}[
  x=1cm,y=1cm,>=Stealth,
  every node/.style={font=\small},
  object/.style={inner sep=3pt},
  map arrow/.style={->,line width=.55pt,shorten <=1pt,shorten >=1pt}
]
\node[object] (M) at (0,2.9) {$B_{\mathbf M_j}$};
\node[object] (Nl) at (-3.2,2.9) {$B_{\mathbf N_{j-1}}$};
\node[object] (Nr) at (3.2,2.9) {$B_{\mathbf N_j}$};
\node[object] (Kl) at (-3.2,1.45) {$B_{\mathcal K_{j-1}}$};
\node[object] (Kr) at (3.2,1.45) {$B_{\mathcal K_j}$};
\node[object] (Il) at (-3.2,0) {$B_{\mathcal I_{j-1}}$};
\node[object] (I) at (0,1.45) {$B_{\mathcal I_j}$};
\node[object] (Ir) at (3.2,0) {$B_{\mathcal I_{j+1}}$};
\node[object] (G) at (0,0) {$(G,*)$};
\draw[map arrow] (M)--node[above=3pt] {$\rho_j$} (Nl);
\draw[map arrow] (M)--node[above=3pt] {$\sigma_j$} (Nr);
\draw[map arrow] (Nl)--node[left=4pt] {$\tau_{j-1}$} (Kl);
\draw[map arrow] (Nr)--node[right=4pt] {$\tau_j$} (Kr);
\draw[map arrow] (Kl)--node[left=4pt] {$a_{j-1}$} (Il);
\draw[map arrow] (Kl)--node[above=3pt] {$b_{j-1}$} (I);
\draw[map arrow] (Kr)--node[above=3pt] {$a_j$} (I);
\draw[map arrow] (Kr)--node[right=4pt] {$b_j$} (Ir);
\draw[map arrow] (Il)--node[above=3pt] {$f_{j-1}$} (G);
\draw[map arrow] (I)--node[right=3pt] {$f_j$} (G);
\draw[map arrow] (Ir)--node[above=3pt] {$f_{j+1}$} (G);
\end{tikzpicture}
\caption{Common refinement at an intermediate representative $f_j$,
for $1\le j<r$.
The two paths from $B_{\mathbf M_j}$ to $B_{\mathcal I_j}$ agree by
\eqref{eq:product-refinement}.
The two lower squares represent the adjacent one-step $C$-homotopies.}
\label{fig:standard-step-refinement}
\end{figure}

At the endpoints,
\[
 f\simeq_C f_0a_0\tau_0,
 \qquad
 f_rb_{r-1}\tau_{r-1}\simeq_C g,
\]
because these representatives are extensions of $f$ and $g$.
Combining these steps gives
\[
 \begin{aligned}
 f&\simeq_C f_0a_0\tau_0\simeq_C f_1b_0\tau_0\\
  &\simeq_C f_1a_1\tau_1\simeq_C\cdots
    \simeq_C f_rb_{r-1}\tau_{r-1}\simeq_C g
    \qquad(r\ge2).
 \end{aligned}
\]
For $r=1$, the extensions $f_0a_0\tau_0$ and $f_1b_0\tau_0$
are joined directly by the step on $B_{\mathbf N_0}$.
Every representative and every common source in this sequence is
standard. Thus inclusion is injective on classes, and (i) gives
surjectivity.

Naturality follows because, for a based digraph map
$f:(G,*)\to(H,*)$, a sphere representative
$f_0:(B_{\mathcal I},\partial B_{\mathcal I})\to(G,*)$,
and every product shrinking $\rho$ with target $B_{\mathcal I}$,
\[
 (f\circ f_0)\rho=f\circ(f_0\rho).\qedhere
\]
\end{proof}

\subsection{Terminal extensions and the colimit}

\begin{definition}\label{def:terminal-shrinking}
For $N\ge n\ge0$, the \emph{terminal shrinking} is
\begin{equation}\label{eq:tail-shrinking}
 p_{N,n}:J_{2N}\longrightarrow J_{2n},
 \qquad
 p_{N,n}(t)=\min(t,2n).
\end{equation}
For a digraph map $\gamma$ on $J_{2n}$,
$(\gamma\circ p_{N,n})(t)=\gamma(t)$ for $t\le2n$ and
$(\gamma\circ p_{N,n})(t)=\gamma(2n)$ for $t\ge2n$.
These shrinkings satisfy
\[
 p_{M,n}=p_{N,n}\circ p_{M,N}
 \qquad(M\ge N\ge n).
\]
\end{definition}

\begin{lemma}\label{lem:pause-normalization}
Let $h:J_{2N}\to J_{2n}$ be a shrinking, with $N\ge n\ge0$.  Then $h$ is homotopic, in the sense of \cref{def:homotopy}, to the terminal shrinking
\[
 p_{N,n}(t)=\min\{t,2n\}.
\]
More precisely, there is a finite sequence of shrinkings
\[
 h=h_0,h_1,\ldots,h_r=p_{N,n}
\]
with fixed endpoint values $h_q(0)=0$ and $h_q(2N)=2n$ for every $q$, such that each adjacent pair is related by a one-step direct or inverse homotopy:
\[
 h_q\toeq h_{q+1}
 \qquad\text{or}\qquad
 h_{q+1}\toeq h_q.
\]
\end{lemma}
\begin{proof}
If $n=0$, both $h$ and $p_{N,0}$ are the constant digraph map to $0$.  Assume $n\ge1$.

For each $0\le j<2n$, surjectivity and monotonicity give a unique transition position
\[
 t_j\in\{0,1,\ldots,2N-1\}
\]
such that
\begin{equation}\label{eq:transition-position}
 h(t_j)=j,
 \qquad
 h(t_j+1)=j+1.
\end{equation}
Hence
\[
 0\le t_0<t_1<\cdots<t_{2n-1}\le2N-1.
\]
The edge between $t_j$ and $t_j+1$ maps to the edge between
$j$ and $j+1$. The alternating orientations therefore give
\begin{equation}\label{eq:transition-parity}
 t_j\equiv j\pmod 2.
\end{equation}
In particular, $t_j-j$ is a nonnegative even integer.

We move the transition positions to $t_j=j$, one at a time.  Suppose that
\[
 t_0=0,\ldots,t_{j-1}=j-1
\]
and $t_j>j$.  By \eqref{eq:transition-parity}, necessarily $t_j\ge j+2$.  Replace $t_j$ by
\[
 t_j'=t_j-2
\]
and leave all other transition positions unchanged.  Then
\[
 t_{j-1}=j-1<t_j'=t_j-2<t_{j+1}
\]
(with the evident omission when $j=0$ or $j=2n-1$), and
\[
 t_j'\equiv j\pmod2.
\]
The shrinking determined by the new list is the map
\[
 h':J_{2N}\to J_{2n},\qquad
 h'(t)=\#\{\,i:0\le i<2n,\ t'_i<t\,\},
 \quad t'_i=t_i\text{ for }i\ne j.
\]
It increases from $0$ to $2n$ by steps of size $1$ at the listed
positions. Since $t'_i\equiv i\pmod2$, each noncollapsed edge has
the orientation of its image.

The two maps differ only at $t_j-1$ and $t_j$, where
\[
 h(t_j-1)=h(t_j)=j,
 \qquad
 h'(t_j-1)=h'(t_j)=j+1.
\]
If $j$ is even, the arrow is $j\to j+1$ and $h\toeq h'$.
If $j$ is odd, it is $j+1\to j$ and $h'\toeq h$.
Thus each move gives a one-step homotopy through shrinkings
fixing both endpoints. After $(t_j-j)/2$ moves the $j$th transition
is in position $j$. Proceeding in order gives $t_j=j$ for all $j$,
which determines $p_{N,n}(t)=\min\{t,2n\}$.
\end{proof}

\begin{example}\label{ex:pauses}
For the shrinking $h_0:J_8\to J_4$, moving the four transitions
in order gives
\[
\begin{aligned}
 h_0&=(0,0,0,1,2,3,4,4,4),\\
 h_1&=(0,1,1,1,2,3,4,4,4),\\
 h_2&=(0,1,2,2,2,3,4,4,4),\\
 h_3&=(0,1,2,3,3,3,4,4,4),\\
 h_4&=(0,1,2,3,4,4,4,4,4)=p_{4,2}.
\end{aligned}
\]
The successive one-step homotopies are
\[
 h_0\toeq h_1\lefteq h_2
 \toeq h_3\lefteq h_4.
\]
Both endpoint values stay fixed throughout.
\end{example}

For $\mathbf n\in\{2,3,\ldots\}^k$, let
$[(B_{\mathbf n},\partial B_{\mathbf n}),(G,*)]$ denote the homotopy
classes of digraph maps of pairs from $(B_{\mathbf n},\partial B_{\mathbf n})$
to $(G,*)$, relative to $\partial B_{\mathbf n}$.
Order these multi-indices coordinatewise. For $\mathbf N\ge\mathbf n$
the product $p_{\mathbf N,\mathbf n}=p_{N_1,n_1}\Box\cdots\Box p_{N_k,n_k}$
defines the transition map $[f]\mapsto[f p_{\mathbf N,\mathbf n}]$.
A relative homotopy pulls back to a relative homotopy because
$p_{\mathbf N,\mathbf n}$ sends boundary to boundary. The composition
law follows from \cref{def:terminal-shrinking}.
In the directed colimit of these sets the equivalence relation is
\begin{equation}\label{eq:colimit-equivalence}
 [f]\sim[g]\quad\Longleftrightarrow\quad
 \exists\,\mathbf M\ge\mathbf n,\mathbf m:\quad
 [f p_{\mathbf M,\mathbf n}]=[g p_{\mathbf M,\mathbf m}].
\end{equation}
Here $f$ and $g$ have sources $B_{\mathbf n}$ and $B_{\mathbf m}$,
respectively, and the equality on the right is taken on $B_{\mathbf M}$.

\begin{proposition}\label{prop:tail-normal-form}
Let $(G,*)$ be a based digraph and let $k\ge1$.
Under the identification of \cref{prop:standard-reduction},
$C$-homotopy is generated by relative homotopies and terminal
extensions in the separate coordinates.  More precisely, let
\[
 f:(B_{\mathbf n},\partial B_{\mathbf n})\longrightarrow(G,*),
 \qquad
 g:(B_{\mathbf m},\partial B_{\mathbf m})\longrightarrow(G,*)
\]
be standard representatives, where
$\mathbf n,\mathbf m\in\{2,3,\ldots\}^k$.  Then
\[
 [f]\sim[g]
 \quad\Longleftrightarrow\quad f\simeq_C g,
\]
where $\sim$ is the colimit relation in \eqref{eq:colimit-equivalence}
and $\simeq_C$ is the relation of \cref{def:sphere-C}.
Consequently there is a natural bijection
\begin{equation}\label{eq:colimit-classes}
 \pC_k(G,*)\cong
 \underset{\mathbf n\in\{2,3,\ldots\}^k}{\operatorname{colim}}
 \bigl[(B_{\mathbf n},\partial B_{\mathbf n}),(G,*)\bigr].
\end{equation}
\end{proposition}
\begin{proof}
Let $\rho=\rho_1\Box\cdots\Box\rho_k:B_{\mathbf N}\to B_{\mathbf n}$
be a product shrinking. By \cref{lem:pause-normalization}, there is a
homotopy through shrinkings from each $\rho_i$ to $p_{N_i,n_i}$, fixing
its endpoint values. Change one coordinate at a time. Adjacent product
maps then have the pointwise comparison in the changing coordinate,
and each map still carries $\partial B_{\mathbf N}$ to $\partial B_{\mathbf n}$.
Postcomposition with $f$ gives
\begin{equation}\label{eq:extension-to-padding}
 f\rho\simeq f p_{\mathbf N,\mathbf n}
 \quad\text{relative to }\partial B_{\mathbf N}.
\end{equation}
Suppose a direct step between standard representatives is witnessed by
$f\rho\toeq g\sigma$ on a common standard box $B_{\mathbf N}$. Then
\[
 f p_{\mathbf N,\mathbf n}\simeq f\rho\toeq g\sigma
   \simeq g p_{\mathbf N,\mathbf m}
 \quad\text{relative to }\partial B_{\mathbf N}.
\]
For an inverse step the middle comparison is $g\sigma\toeq f\rho$.
If $f\simeq_C g$, use \cref{prop:standard-reduction} to choose a finite
sequence of such steps on standard boxes. Take $\mathbf M$ to be the
coordinatewise maximum of the common-source indices of these steps.
Pull the homotopy on each $B_{\mathbf N}$ back along
$p_{\mathbf M,\mathbf N}$. If an intermediate representative $f_j$
has source $B_{\mathbf n_j}$, its two extensions to $B_{\mathbf M}$
both equal $f_j p_{\mathbf M,\mathbf n_j}$ by
$p_{\mathbf N,\mathbf n_j}\circ p_{\mathbf M,\mathbf N}
=p_{\mathbf M,\mathbf n_j}$. The composite
homotopy proves $[f]\sim[g]$ by \eqref{eq:colimit-equivalence}.
For a sequence of length zero this conclusion holds with $f=g$.

Conversely, a relative homotopy on $B_{\mathbf M}$ as in
\eqref{eq:colimit-equivalence} gives
\[
 f\simeq_C f p_{\mathbf M,\mathbf n}\simeq_C
 g p_{\mathbf M,\mathbf m}\simeq_C g
\]
by \eqref{eq:extension-C}. Every sphere class has a standard
representative, and postcomposition by a based map commutes with
terminal extension. This proves the stated natural bijection.
\end{proof}

\begin{remark}\label{rem:sphere-condition-standardization}
Admitting all products of finite oriented lines of positive length
gives the same $C$-homotopy classes by the proofs of
\cref{prop:standard-reduction,prop:tail-normal-form}.
The condition in \cref{def:digraph-sphere} ensures that each admitted
quotient has the integral GLMY path homology of a sphere of the given
dimension.
\end{remark}

\section{Sphere homotopy groups}\label{sec:groups}

We use the description of $C$-homotopy in \cref{prop:tail-normal-form}
to show that concatenation in the first coordinate gives a product on
sphere classes. This product is a group law, and it is commutative
for $k\ge2$.

In the formulas below, $B_{(n_2,\ldots,n_k)}$ is the standard box in
coordinates $2,\ldots,k$.  For $k=1$, the empty product is interpreted
as a one-vertex digraph with empty coordinate boundary.

\begin{definition}\label{def:concat}
Let $(G,*)$ be a based digraph and $k\ge1$.
After terminal extension to common lengths $2n_i$ in coordinates
$i=2,\ldots,k$, let the two standard representatives be the digraph maps of pairs
\[
 f:\bigl(J_{2m}\Box B_{(n_2,\ldots,n_k)},
         \partial(J_{2m}\Box B_{(n_2,\ldots,n_k)})\bigr)\to(G,*),
\]
and
\[
 g:\bigl(J_{2n}\Box B_{(n_2,\ldots,n_k)},
         \partial(J_{2n}\Box B_{(n_2,\ldots,n_k)})\bigr)\to(G,*).
\]
Their \emph{first-coordinate concatenation} is the digraph map of pairs
\[
 f\#_1 g:
 \bigl(J_{2(m+n)}\Box B_{(n_2,\ldots,n_k)},
       \partial(J_{2(m+n)}\Box B_{(n_2,\ldots,n_k)})\bigr)\to(G,*)
\]
defined by
\begin{equation}\label{eq:concat}
 (f\#_1g)(t,x)=
 \begin{cases}
 f(t,x),&0\le t\le2m,\\
 g(t-2m,x),&2m\le t\le2m+2n.
 \end{cases}
\end{equation}
The formulas agree at $t=2m$ since $f(2m,x)=g(0,x)=*$.
Each arrow lies in one of the two subboxes, and the even shift by
$2m$ preserves its orientation. Both maps take their boundaries to $*$,
so $f\#_1g$ is a digraph map of pairs.
\end{definition}

\begin{proposition}\label{prop:concat}
Concatenation induces a well-defined product on $\pC_k(G,*)$, independent of the terminal extensions in coordinates $2,\ldots,k$ and of the chosen representatives.
\end{proposition}
\begin{proof}
By \cref{prop:tail-normal-form}, it suffices to check relative
homotopies and terminal extensions.  Take $f,g$ as in \cref{def:concat}.

Suppose $f_0\toeq f_1$ on $J_{2m}\Box B_{(n_2,\ldots,n_k)}$.
On the first subbox of \eqref{eq:concat}, the concatenations
$f_0\#_1g$ and $f_1\#_1g$ have this pointwise comparison. On the
second subbox they agree. Comparing the second factor in the same way gives
\[
 f_0\toeq f_1\ \Longrightarrow\ f_0\#_1g\toeq f_1\#_1g,
 \qquad
 g_0\toeq g_1\ \Longrightarrow\ f\#_1g_0\toeq f\#_1g_1.
\]
Both comparisons fix the full boundary. Reversing the arrows gives
the inverse steps.

For $k\ge2$ and $(N_2,\ldots,N_k)\ge(n_2,\ldots,n_k)$, coordinatewise terminal
extension gives
\[
\begin{aligned}
 &\bigl(f\circ(\id_{J_{2m}}\Box p_{(N_2,\ldots,N_k),(n_2,\ldots,n_k)})\bigr)\#_1
  \bigl(g\circ(\id_{J_{2n}}\Box p_{(N_2,\ldots,N_k),(n_2,\ldots,n_k)})\bigr)\\
 &\qquad=(f\#_1g)\circ
  (\id_{J_{2(m+n)}}\Box p_{(N_2,\ldots,N_k),(n_2,\ldots,n_k)})\\
 &\qquad\simeq_C f\#_1g.
\end{aligned}
\]
Any two choices of lengths in coordinates $2,\ldots,k$ admit a common
terminal extension, so this identity proves independence of those choices.

For an integer $d\ge0$, define
\[
 r:J_{2(m+n+d)}\longrightarrow J_{2(m+n)},\qquad
 r(t)=
 \begin{cases}
  t,&0\le t\le2m,\\
  2m,&2m\le t\le2m+2d,\\
  t-2d,&2m+2d\le t\le2(m+n+d).
 \end{cases}
\]
The map $r$ is nondecreasing and onto, and the even shift by $2d$
preserves the orientation of each shifted edge. Thus $r$ is a
shrinking. By \eqref{eq:concat} and \eqref{eq:extension-C},
\[
\begin{aligned}
 &\bigl(f\circ(p_{m+d,m}\Box\id_{B_{(n_2,\ldots,n_k)}})\bigr)\#_1g\\
 &\qquad=(f\#_1g)\circ(r\Box\id_{B_{(n_2,\ldots,n_k)}})
 \simeq_C f\#_1g.
\end{aligned}
\]
For terminal extension of the second factor,
\[
\begin{aligned}
 &f\#_1\bigl(g\circ(p_{n+d,n}\Box\id_{B_{(n_2,\ldots,n_k)}})\bigr)\\
 &=(f\#_1g)\circ(p_{m+n+d,m+n}\Box\id_{B_{(n_2,\ldots,n_k)}})
 \simeq_C f\#_1g.
\end{aligned}
\]
These identities show that concatenation respects relative
homotopies and terminal extensions.  By \cref{prop:tail-normal-form},
the product $[f]\cdot[g]=[f\#_1g]$ is therefore well-defined.
\end{proof}

\begin{theorem}\label{thm:sphere-groups}
The product induced by $\#_1$ makes $\pC_k(G,*)$ a group.  Its identity is
represented by any constant digraph map of pairs on a $k$-sphere
presentation.  If
\[
 f:\bigl(J_{2m}\Box B_{(n_2,\ldots,n_k)},
          \partial(J_{2m}\Box B_{(n_2,\ldots,n_k)})\bigr)\longrightarrow(G,*)
\]
is a standard representative, then the inverse of $[f]$ is represented by
the digraph map of pairs
\[
 f^{-}:\bigl(J_{2m}\Box B_{(n_2,\ldots,n_k)},
          \partial(J_{2m}\Box B_{(n_2,\ldots,n_k)})\bigr)\longrightarrow(G,*),
\]
defined by
\[
 f^{-}(t,x)=f(2m-t,x).
\]
\end{theorem}
\begin{proof}
Let $[f],[g],[h]\in\pC_k(G,*)$.  After terminal extension in coordinates
$2,\ldots,k$, choose representatives
\[
\begin{aligned}
 f&:J_{2m}\Box B_{(N_2,\ldots,N_k)}\longrightarrow G,\\
 g&:J_{2n}\Box B_{(N_2,\ldots,N_k)}\longrightarrow G,\\
 h&:J_{2r}\Box B_{(N_2,\ldots,N_k)}\longrightarrow G,
\end{aligned}
\]
all constant on the full coordinate boundary.  Both
$(f\#_1g)\#_1h$ and $f\#_1(g\#_1h)$ are digraph maps on
$J_{2(m+n+r)}\Box B_{(N_2,\ldots,N_k)}$, and both are given by
\begin{equation}\label{eq:assoc-piecewise}
 (t,x)\longmapsto
 \begin{cases}
  f(t,x),&0\le t\le2m,\\
  g(t-2m,x),&2m\le t\le2m+2n,\\
  h(t-2m-2n,x),&2m+2n\le t\le2m+2n+2r.
 \end{cases}
\end{equation}
Hence
\[
 (f\#_1g)\#_1h=f\#_1(g\#_1h)
\]
strictly after the common terminal extensions, and therefore
\[
 ([f]\cdot[g])\cdot[h]=[f]\cdot([g]\cdot[h]).
\]

All constant sphere representatives have the same class. Indeed,
\cref{prop:standard-reduction} extends each to a constant standard
representative, and two such representatives agree after terminal
extension to the coordinatewise maximum of their indices.
Let $e_s:J_{2s}\Box B_{(N_2,\ldots,N_k)}\to G$ be the constant digraph map
with value $*$, where $s\ge2$, and write $[e]=[e_s]$.  Then
\[
\begin{aligned}
 (f\#_1e_s)(t,x)&=f(\min(t,2m),x),\\
 (e_s\#_1f)(t,x)&=f(\max(0,t-2s),x).
\end{aligned}
\]
Both substitutions $J_{2(m+s)}\to J_{2m}$ are shrinkings, since
each noncollapsed edge retains its position or shifts by $2s$.
Thus both concatenations are extensions of $f$ along product shrinkings, and
\[
 [f]\cdot[e]=[f]=[e]\cdot[f].
\]

The reflection
\[
 \iota_m:J_{2m}\longrightarrow J_{2m},\qquad
 \iota_m(t)=2m-t,
\]
preserves parity and adjacency, so it is a digraph isomorphism of
$J_{2m}$ and exchanges the two endpoints. Hence
\[
 f^{-}=f\circ(\iota_m\Box\id_{B_{(N_2,\ldots,N_k)}})
\]
is a digraph map of pairs.

To contract $f\#_1f^{-}$, we successively remove the two copies of
the last edge traversed in the first coordinate, as in \cite{LiMY}.
For $q=2m,2m-1,\ldots,0$, define
\[
 H_q:\bigl(J_{4m}\Box B_{(N_2,\ldots,N_k)},
       \partial(J_{4m}\Box B_{(N_2,\ldots,N_k)})\bigr)\longrightarrow(G,*)
\]
by
\[
 H_q(t,x)=
 \begin{cases}
 f(t,x),&0\le t\le q,\\
 f(q,x),&q\le t\le4m-q,\\
 f(4m-t,x),&4m-q\le t\le4m.
 \end{cases}
\]
Then
\[
 H_{2m}=f\#_1f^{-},\qquad H_0(t,x)=f(0,x)=*.
\]
On the three subboxes, $H_q$ follows $f$, stays at the slice
$f(q,\cdot)$, and follows $f$ in reverse. The formulas agree where
the subboxes meet and take the full boundary to $*$. The first and
last formulas use the inclusion and the parity-preserving reflection
$t\mapsto4m-t$, and the middle formula collapses all first-coordinate
edges. Thus every $H_q$ is a digraph map of pairs.
For $q\le t\le4m-q$, the values of $H_q$ and $H_{q-1}$ are $f(q,x)$
and $f(q-1,x)$, respectively. The maps agree for all other $t$.
Since $f$ is a digraph map, the alternating arrows of
$J_{2m}$ give
\[
 \begin{cases}
 H_{q-1}\toeq H_q,&q\text{ odd},\\
 H_q\toeq H_{q-1},&q\text{ even},
 \end{cases}
 \qquad 1\le q\le2m.
\]
Thus each adjacent pair $H_q,H_{q-1}$ is related by a one-step direct or
inverse relative homotopy.  Hence
\[
 [f]\cdot[f^{-}]=[H_{2m}]=[H_0]=[e].
\]
Applying the same construction to $f^{-}$ gives
\[
 [f^{-}]\cdot[f]=[e].
\]
Thus every class has a two-sided inverse, and $\pC_k(G,*)$ is a group.
\end{proof}

\begin{proposition}\label{prop:abelian}
For $k\ge2$, the group $\pC_k(G,*)$ is abelian.
\end{proposition}
\begin{proof}
We use the Eckmann--Hilton interchange argument \cite{EckmannHilton}.
For $k\ge2$, concatenation in the first and second coordinates induces
products $\cdot_1$ and $\cdot_2$ on $\pC_k(G,*)$, respectively.  The proofs
of \cref{prop:concat,thm:sphere-groups} apply in both coordinates, so the
products are well defined and have the same identity element $[e]$.

Take $[f],[g],[h],[\ell]\in\pC_k(G,*)$.  After terminal extension, choose
representative digraph maps of pairs on one common standard alternating box,
\[
 f,g,h,\ell:
 (B_{\mathbf M},\partial B_{\mathbf M})
 \longrightarrow(G,*),
 \qquad\mathbf M=(M_1,\ldots,M_k).
\]
Both iterated concatenations
\[
 (f\#_1g)\#_2(h\#_1\ell),
 \qquad
 (f\#_2h)\#_1(g\#_2\ell)
\]
are digraph maps on
\[
 \begin{cases}
 J_{4M_1}\Box J_{4M_2},&k=2,\\
 J_{4M_1}\Box J_{4M_2}\Box J_{2M_3}\Box\cdots\Box J_{2M_k},&k\ge3.
 \end{cases}
\]
At every vertex $(x_1,\ldots,x_k)$, both digraph maps have the value
\begin{equation}\label{eq:interchange-blocks}
 \begin{cases}
 f(x_1,\ldots,x_k),&0\le x_1\le2M_1,\ 0\le x_2\le2M_2,\\
 g(x_1-2M_1,x_2,\ldots,x_k),&2M_1\le x_1\le4M_1,\ 0\le x_2\le2M_2,\\
 h(x_1,x_2-2M_2,x_3,\ldots,x_k),&0\le x_1\le2M_1,\ 2M_2\le x_2\le4M_2,\\
 \ell(x_1-2M_1,x_2-2M_2,x_3,\ldots,x_k),&2M_1\le x_1\le4M_1,\ 2M_2\le x_2\le4M_2.
 \end{cases}
\end{equation}
The two maps thus agree pointwise, giving the interchange identity
on classes
\begin{equation}\label{eq:interchange-law}
 ([f]\cdot_1[g])\cdot_2([h]\cdot_1[\ell])
 =
 ([f]\cdot_2[h])\cdot_1([g]\cdot_2[\ell]).
\end{equation}

Taking $[g]=[h]=[e]$ in \eqref{eq:interchange-law} gives
\begin{equation}\label{eq:two-products-equal}
 [f]\cdot_2[\ell]=[f]\cdot_1[\ell].
\end{equation}
Taking $[f]=[\ell]=[e]$ in \eqref{eq:interchange-law} gives
\[
 [g]\cdot_2[h]=[h]\cdot_1[g].
\]
Using \eqref{eq:two-products-equal},
\[
 [g]\cdot_1[h]=[g]\cdot_2[h]=[h]\cdot_1[g].
\]
Hence the group operation is commutative.
\end{proof}

We use multiplicative notation in formulas valid for all $k\ge1$.
For $k\ge2$, we also write $[f]+[g]=[f\#_1g]$ and $0=[e]$.

\begin{proposition}\label{prop:functoriality}
Let $k\ge1$ and let $f:(G,*)\to(H,*)$ be a based digraph map.
Postcomposition defines a group homomorphism $f_*:\pC_k(G,*)\to\pC_k(H,*)$.
If the based digraph maps $f,g:(G,*)\to(H,*)$ are homotopic
relative to the basepoint, then $f_*=g_*$.  In particular, based homotopy equivalences induce
isomorphisms.
\end{proposition}
\begin{proof}
Define $f_*([f_0])=[f\circ f_0]$. Postcomposition preserves the
pointwise comparisons witnessing each $C$-homotopy step, since
$f_0\rho\toeq g_0\sigma$ implies
$(f\circ f_0)\rho\toeq(f\circ g_0)\sigma$ on the same source.
It also preserves the product, since on standard representatives
\[
 f\circ(f_0\#_1 g_0)=(f\circ f_0)\#_1(f\circ g_0).
\]
If $F:G\Box I_m\to H$ is a based homotopy from $f$ to $g$,
then $F\circ(f_0\Box\id_{I_m})$ is a homotopy from $f\circ f_0$
to $g\circ f_0$ relative to the full boundary of the source of $f_0$.
Thus $f_*=g_*$. Identity maps and compositions act by postcomposition,
so based homotopy inverses induce inverse isomorphisms.
\end{proof}

\begin{remark}\label{rem:cubical-comparison}
The colimit in \eqref{eq:colimit-classes} gives the cubical homotopy
classes $A_k(G)$ of
\cite[Definition~1.3.1 and Corollary~1.3.7, equation~(1.3.6)]{EldridgeEtAl}.
For the comparison with that paper, we adjoin a reflexive arrow
at each vertex of our digraphs, as required by their convention.
The standard intervals in that paper have the same
alternating orientations as our standard lines at even lengths.
Their right shrinkings send $t$ to $\min\{t,m\}$ when the target
interval has length $m$. In that right-shrinking sequence we can restrict to lengths
$2N$ with $N\ge2$. In \eqref{eq:colimit-classes} we can restrict to
the diagonal indices $(N,\ldots,N)$, since they are cofinal in the
coordinatewise order. The relative homotopy classes at these stages
and their transition maps are then the same in both colimits.

Choose representatives $f$ and $g$ on a common standard cube
$J_{2N}^{\Box k}$ by terminal extension. The group law in
\cref{thm:sphere-groups} and the elementary group law cited in
\cite[Remark~1.3.3]{EldridgeEtAl} both represent the product by
$f\#_1g$. Its terminal extension in coordinates $2,\ldots,k$
to $J_{4N}^{\Box k}$ represents the same product in the diagonal colimit.
Since the identification
commutes with postcomposition, it gives a natural group isomorphism
\[
 \pC_k(G,*)\xrightarrow{\cong}A_k(G),\qquad k\ge1.
\]
\end{remark}

\begin{remark}\label{rem:symmetric-targets}
Let $(X,*)$ be a based simple graph, and let $\overleftrightarrow X$
be the digraph obtained by giving each edge of $X$ both orientations.
There is a natural isomorphism
\[
 \pC_k(\overleftrightarrow X,*)\cong A_k(X,*)\qquad(k\ge1),
\]
where $A_k$ is defined by based maps on the lattice $\Z^k$
\cite[Definition~1.18]{CarranzaKapulkin}.
On a fixed box, a map to $\overleftrightarrow X$ is precisely a
graph map from the underlying grid to $X$ allowing edge collapses,
since each edge of $X$ has both orientations in $\overleftrightarrow X$. The same
observation applies to relative homotopies. In the
colimit of \cref{prop:tail-normal-form}, take the cofinal cubes
$J_{2R}^{\Box k}$ with even $R\ge2$. Regard each underlying grid
as $[-R,R]^k\cap\Z^k$. For even $S\ge R$, the coordinate map
\[
 J_{2S}\longrightarrow J_{2R},\qquad
 t\longmapsto\min\{2R,\max\{0,t-(S-R)\}\}
\]
is a shrinking because $S-R$ is even. By \cref{lem:pause-normalization}
it gives the same transition on relative classes as terminal
extension. In centered coordinates it pads the cube equally on both
sides. Extension by $*$ to $\Z^k$ now commutes with the transitions.
Conversely, a based lattice map and all the time slices of a based
lattice homotopy are constant at $*$ off a common finite cube.
Restriction to a larger centered cube therefore gives the inverse on
classes. For representatives $f,g$ on $J_{2R}^{\Box k}$, let
$\widetilde f,\widetilde g$ denote their extensions by $*$ in the
centered coordinates, and use the same convention for $f\#_1g$.
The lattice product $\widetilde f\cdot_1\widetilde g$ in
\cite[Definition~1.18]{CarranzaKapulkin} agrees with $\widetilde f$
for $t\le R$ and with $\widetilde g(t-2R,x)$ for $t>R$, where $t$
is the first coordinate and $x$ denotes the others. Hence
\[
 \widetilde{f\#_1g}(t,x)=(\widetilde f\cdot_1\widetilde g)(t+R,x)
 \qquad ((t,x)\in\Z^k).
\]
Precomposing by the successive shifts $t\mapsto t+j$, $0\le j\le R$,
gives a based homotopy: consecutive inputs are adjacent, and all the
translated maps equal $*$ outside a common finite cube. The
identification preserves multiplication and commutes with based graph maps.
\end{remark}

\begin{proposition}\label{prop:pi1}
Taking a one-dimensional sphere representative to its underlying based
loop induces a natural group isomorphism
\[
 \alpha_1^G:\pC_1(G,*)\xrightarrow{\cong}\pi_1(G,*),
 \qquad \alpha_1^G([f])=[f]_{\mathrm{GLMY}}.
\]
\end{proposition}
\begin{proof}
A direct sphere $C$-homotopy from $f$ to $g$ with shrinkings
$\rho,\sigma$ gives
\[
 f\toeq f\circ\rho\toeq g\circ\sigma\toeq g
 \qquad\text{in }L(G,*)
\]
by \eqref{eq:classical-arrow}, so $\alpha_1^G$ is well defined.

Every based loop $\gamma$ on $I_m$, including the length-zero loop,
has a standard extension $\gamma\circ r$ along a shrinking
$r:J_{2N}\to I_m$ with $N\ge\max\{m,2\}$
(\cref{lem:standard-cover,thm:admissibility}).
Any two such extensions agree after common refinement
(\cref{lem:common-refinement}). To compare different based loops,
let $\gamma:I_M\to G$ and $\delta:I_m\to G$ with $M\ge m$.
For a GLMY arrow between $\gamma$ and $\delta$, let $h:I_M\to I_m$
be the shrinking in \eqref{eq:classical-arrow}. If $M=m$, name
the source loop $\gamma$. Choose a shrinking $s:J_{2N}\to I_M$
with $N\ge\max\{M,2\}$. Composing the arrow comparison with $s$ gives
\[
 \begin{cases}
  \gamma\circ s\toeq\delta\circ h\circ s,&\gamma\to\delta,\\
  \delta\circ h\circ s\toeq\gamma\circ s,&\delta\to\gamma .
 \end{cases}
\]
Thus the two standard extensions are $C$-homotopic.
Thus $[\gamma]_{\mathrm{GLMY}}\mapsto[\gamma\circ r]$
is well defined and inverts $\alpha_1^G$, since extensions preserve
both classes. On standard representatives, $\#_1$ is loop
concatenation, and $\alpha_1^G$ commutes with postcomposition.
\end{proof}

\section{The Hurewicz homomorphism}\label{sec:top-class}

We construct an oriented top cycle on each standard digraph sphere.
A standard representative sends this cycle to a cycle in $G$.
We show that the resulting class in integral GLMY path homology
depends only on the sphere homotopy class and is additive under
concatenation.

\begin{definition}\label{def:signed-line-chain}
For a finite oriented line $I_m$, define
\begin{equation}\label{eq:line-chain}
 c_{I_m}
 =\sum_{i\to i+1}e_{i,i+1}
 -\sum_{i+1\to i}e_{i+1,i}.
\end{equation}
Then
\[
 \partial c_{I_m}=e_m-e_0.
\]
An arrow has coefficient $+1$ if it points in the order $0,\ldots,m$
and $-1$ otherwise. In particular,
\begin{equation}\label{eq:standard-chain}
 c_{J_{2n}}=\sum_{r=0}^{n-1}
 (e_{2r,2r+1}-e_{2r+2,2r+1}).
\end{equation}
\end{definition}

\begin{definition}\label{def:staircase-product}
Let $X,Y$ be digraphs and $p,q\ge0$.  Let $e_{x_0\cdots x_p}$ and
$e_{y_0\cdots y_q}$ be regular paths in $X$ and $Y$, respectively.  A
$(p,q)$-shuffle is a word $\sigma$ with $p$ copies of the letter $X$
and $q$ copies of the letter $Y$. Write $\operatorname{Sh}(p,q)$ for
the set of such words. Starting at $(i_0,j_0)=(0,0)$, let
$i_\ell$ and $j_\ell$ count the letters $X$ and $Y$ among the first
$\ell$ letters of $\sigma$. The resulting path in $X\Box Y$ has vertex
sequence
\[
 (x_{i_0},y_{j_0}),\ (x_{i_1},y_{j_1}),\ \ldots,\
 (x_{i_{p+q}},y_{j_{p+q}}).
\]
The path is called \emph{step-like} because each step changes
exactly one coordinate.
Let $\operatorname{inv}(\sigma)$ be the number of
pairs in which a $Y$ occurs before an $X$.  Define
\[
e_{x_0\cdots x_p}\times e_{y_0\cdots y_q}
 =\sum_{\sigma\in\operatorname{Sh}(p,q)}
 (-1)^{\operatorname{inv}(\sigma)}
 e_{(x_{i_0},y_{j_0})\cdots(x_{i_{p+q}},y_{j_{p+q}})},
\]
and extend bilinearly to
\[
 \times:\mathcal R_p(X;\Z)\times\mathcal R_q(Y;\Z)
 \longrightarrow\mathcal R_{p+q}(X\Box Y;\Z).
\]
This is the \emph{cross product} of paths. For more than two factors
we apply it from left to right, identifying both bracketings of
a box product with the same product of vertex sets. For homogeneous
chains $u$ and $v$, with
$|u|$ denoting the degree of $u$, the Leibniz rule
\cite[Definition~4.2 and Proposition~4.4]{GMY} is
\[
 \partial(u\times v)=\partial u\times v+(-1)^{|u|}u\times\partial v.
\]
For two arrows the shuffle formula gives
\begin{equation}\label{eq:two-shuffles}
 e_{x_0x_1}\times e_{y_0y_1}
 =e_{(x_0,y_0)(x_1,y_0)(x_1,y_1)}
 -e_{(x_0,y_0)(x_0,y_1)(x_1,y_1)}.
\end{equation}
\end{definition}

For a standard alternating box
$B_{\nn}=J_{2n_1}\Box\cdots\Box J_{2n_k}$, write
\[
 \begin{aligned}
 C_{\nn}&=c_{J_{2n_1}}\times\cdots\times c_{J_{2n_k}},\\
 z_{\nn}&=(q_{\nn})_\#C_{\nn}.
 \end{aligned}
\]

\begin{lemma}\label{lem:standard-top-chain}
These chains satisfy
\[
 C_{\nn}\in\Omega_k(B_{\nn};\Z),\qquad
 z_{\nn}\in\Omega_k(\Sn;\Z),\qquad
 \partial z_{\nn}=0.
\]
\end{lemma}
\begin{proof}
Every step-like path occurring in the cross product of allowed paths is allowed, so
\[
 C_{\nn}\in\mathcal A_k(B_{\nn};\Z).
\]
Since $\partial c_{J_{2n_i}}=e_{2n_i}-e_0$, the Leibniz rule in
\cref{def:staircase-product} gives
\[
 \begin{aligned}
 \partial C_{\nn}
 &=\sum_{i=1}^k(-1)^{i-1}
 c_{J_{2n_1}}\times\cdots\times(e_{2n_i}-e_0)
       \times\cdots\times c_{J_{2n_k}}\\
 &\in\mathcal A_{k-1}(B_{\nn};\Z).
 \end{aligned}
\]
Thus $C_{\nn}\in\Omega_k(B_{\nn};\Z)$ by \cref{def:path-complex}.
The quotient digraph map $q_{\nn}$ induces a chain map, so
$z_{\nn}\in\Omega_k(\Sn;\Z)$.
The boundary of $C_{\nn}$ is supported on the coordinate faces.
When $k>1$, the boundary collapse sends every elementary path in
$\partial C_{\nn}$ to a constant path of degree $k-1$, which
vanishes in the regular path module.
When $k=1$, the two endpoint terms have the same image and opposite
signs.  Therefore
\[
 \partial z_{\nn}=(q_{\nn})_\#\partial C_{\nn}=0.\qedhere
\]
\end{proof}

\begin{lemma}\label{lem:standard-generator}
For every standard digraph sphere $(\Sn,*)$, the cycle
\[
 z_{\nn}=(q_{\nn})_\#C_{\nn}
\]
generates $H_k(\Sn;\Z)\cong\Z$.
\end{lemma}
\begin{proof}
An allowed path in the interior changes each coordinate at most once,
since every arrow of an alternating line ends at a sink. No arrow
enters $*$. If an allowed path starts at $*$, its first arrow reaches
a vertex whose $i$th coordinate is $1$ or $2n_i-1$ for some $i$.
The value $1$ or $2n_i-1$ is a sink of $J_{2n_i}$, so coordinate
$i$ remains fixed after the first arrow. Hence at most $k-1$ arrows
follow the first.
Thus $\mathcal A_{k+1}(\Sn;\Z)=0$, and the coefficient of any given
elementary $k$-path defines a homomorphism $H_k(\Sn;\Z)\to\Z$.

Consider the path from $(2,\ldots,2)$ to $(3,\ldots,3)$ which changes
the coordinates in the order $1,\ldots,k$. All its vertices are interior
because $n_i\ge2$, so it has a unique lift to the box. The factor
$e_{2,3}$ has coefficient $+1$ in each signed line chain, and the
displayed order has shuffle sign $+1$. Hence the coefficient of this
path in $z_{\nn}$ is $+1$.
Since $H_k(\Sn;\Z)\cong\Z$ by \cref{thm:admissibility}, the
coefficient homomorphism can take the value $1$ on $[z_{\nn}]$ only
if $[z_{\nn}]$ generates $H_k(\Sn;\Z)$.
\end{proof}

By \cref{def:sphere-representative,lem:standard-top-chain},
each standard representative $f$ gives a cycle
$f_\#C_{\nn}=\bar f_\#z_{\nn}$ in $\Omega_k(G;\Z)$.
We compare these cycles under changes of representative.

\begin{lemma}\label{lem:top-shrinking}
For a product shrinking $\rho:B_{\mathbf N}\to B_{\nn}$ of standard
boxes, the endpoints in each coordinate map to the corresponding
endpoints. Hence $\rho$ induces a based digraph map
$\bar\rho:(S^k_{2N_1\times\cdots\times2N_k},*)\to(\Sn,*)$ with
\[
  \bar\rho\circ q_{\mathbf N}=q_{\nn}\circ\rho.
\]
Then
\[
  \rho_\#C_{\mathbf N}=C_{\nn},\qquad
  \bar\rho_\#z_{\mathbf N}=z_{\nn}.
\]
\end{lemma}
\begin{proof}
Write $\rho=\rho_1\Box\cdots\Box\rho_k$, where
$\rho_i:J_{2N_i}\to J_{2n_i}$ are shrinkings. By
\eqref{eq:transition-position} and \eqref{eq:transition-parity},
for each $0\le j<2n_i$ there is a unique $t$ with
$\rho_i(t)=j$ and $\rho_i(t+1)=j+1$, and $t\equiv j\pmod2$.
Thus the signed source edge at this transition contributes the
coefficient of the signed target edge at $j$. Using $e_{aa}=0$ in
the regular path module gives
\[
 \begin{aligned}
 (\rho_i)_\#c_{J_{2N_i}}
 &=\sum_{r=0}^{N_i-1}
   \bigl(e_{\rho_i(2r),\rho_i(2r+1)}
        -e_{\rho_i(2r+2),\rho_i(2r+1)}\bigr)\\
 &=\sum_{r=0}^{n_i-1}
   \bigl(e_{2r,2r+1}-e_{2r+2,2r+1}\bigr)
 =c_{J_{2n_i}}.
 \end{aligned}
\]
For coordinate arrows $a_i\to b_i$ in $J_{2N_i}$, applying $\rho$ to
each step-like path in the cross product gives
\[
 \rho_\#(e_{a_1b_1}\times\cdots\times e_{a_kb_k})
 =\begin{cases}
  0,&\rho_i(a_i)=\rho_i(b_i)\text{ for some }i,\\
  e_{\rho_1(a_1),\rho_1(b_1)}\times\cdots\times
  e_{\rho_k(a_k),\rho_k(b_k)},&\text{otherwise}.
 \end{cases}
\]
If an arrow collapses, every image path has consecutive equal vertices.
Otherwise the order of the coordinate steps and its sign are preserved.
Multilinearity gives
\[
 \begin{aligned}
 \rho_\#C_{\mathbf N}
 &=\rho_\#(c_{J_{2N_1}}\times\cdots\times c_{J_{2N_k}})\\
 &=(\rho_1)_\#c_{J_{2N_1}}\times\cdots\times
   (\rho_k)_\#c_{J_{2N_k}}\\
 &=c_{J_{2n_1}}\times\cdots\times c_{J_{2n_k}}
 =C_{\nn}.
 \end{aligned}
\]
Finally, $\bar\rho\circ q_{\mathbf N}=q_{\nn}\circ\rho$ implies
\begin{align*}
 \bar\rho_\#z_{\mathbf N}
 &=\bar\rho_\#(q_{\mathbf N})_\#C_{\mathbf N}\\
 &=(q_{\nn})_\#\rho_\#C_{\mathbf N}\\
 &=(q_{\nn})_\#C_{\nn}=z_{\nn}.\qedhere
\end{align*}
\end{proof}

\begin{lemma}\label{lem:relative-prism}
Let $f_0,f_1:(B_{\nn},\partial B_{\nn})\to(G,*)$ be digraph maps of pairs
that are homotopic relative to $\partial B_{\nn}$.
There is a chain $P\in\Omega_{k+1}(G;\Z)$ such that
\begin{equation}\label{eq:prism-certificate}
 \partial P=(f_1)_\#C_{\nn}-(f_0)_\#C_{\nn}.
\end{equation}
\end{lemma}
\begin{proof}
Let the relative homotopy be represented by a digraph map of pairs
\[
 H:\bigl(B_{\nn}\Box I_m,\,\partial B_{\nn}\Box I_m\bigr)
 \longrightarrow(G,*),\qquad
 H(x,0)=f_0(x),\quad H(x,m)=f_1(x).
\]
Every term in $C_{\nn}\times c_{I_m}$ is allowed. Its boundary is
\[
 (\partial C_{\nn})\times c_{I_m}
   +(-1)^k C_{\nn}\times(e_m-e_0).
\]
Both summands are allowed by \cref{lem:standard-top-chain} and the
shuffle formula, so $C_{\nn}\times c_{I_m}$ is $\partial$-invariant.
Define the prism chain by
\[
 P=(-1)^k H_\#(C_{\nn}\times c_{I_m}).
\]
The endpoint terms are the images of $C_{\nn}\times e_0$ and
$C_{\nn}\times e_m$ under $H_\#$, with the endpoint signs given by
$e_m-e_0$. Thus
\[
 \partial P=(-1)^kH_\#\bigl((\partial C_{\nn})\times c_{I_m}\bigr)
        +(f_1)_\#C_{\nn}-(f_0)_\#C_{\nn}.
\]
Each elementary path in $(\partial C_{\nn})\times c_{I_m}$ lies in
$\partial B_{\nn}\Box I_m$, so $H$ sends all its vertices to $*$.
A constant path of degree $k\ge1$ vanishes in the regular path
module. This proves
\eqref{eq:prism-certificate}.
\end{proof}

\begin{samepage}
\begin{theorem}[Hurewicz homomorphism]\label{thm:top-map}
For every based digraph $(G,*)$ and every $k\ge1$, the assignment
\begin{equation}\label{eq:top-map}
 h_k^C:\pC_k(G,*)\longrightarrow H_k(G;\Z),\qquad
 [f]\longmapsto H_k(\bar f)([z_{\nn}])=[f_\#C_{\nn}]
\end{equation}
is a well-defined natural group homomorphism, where
$f:(B_{\nn},\partial B_{\nn})\to(G,*)$ is a standard representative.
\end{theorem}
\end{samepage}
\begin{proof}
Every class has a standard representative by \cref{prop:standard-reduction}.
Suppose $f$ on $B_{\mathbf n}$ and $g$ on $B_{\mathbf m}$ represent
the same class. By \cref{prop:tail-normal-form}, there is a common
standard box $B_{\mathbf N}$ on which $f p_{\mathbf N,\mathbf n}$
and $g p_{\mathbf N,\mathbf m}$ are homotopic relative to the boundary.
Apply \cref{lem:relative-prism} to that homotopy. Its prism chain
$P\in\Omega_{k+1}(G;\Z)$ satisfies
\begin{equation}\label{eq:top-invariance}
 \begin{aligned}
 \partial P
 &=(g p_{\mathbf N,\mathbf m})_\#C_{\mathbf N}
   -(f p_{\mathbf N,\mathbf n})_\#C_{\mathbf N}\\
 &=g_\#C_{\mathbf m}-f_\#C_{\mathbf n},
 \end{aligned}
\end{equation}
where the second equality follows from \cref{lem:top-shrinking}.
Thus \eqref{eq:top-map} is well defined.

To compute $h_k^C$ on a product, use terminal extensions in
coordinates $2,\ldots,k$ to choose standard representatives
\[
 \begin{aligned}
  f&:\bigl(J_{2r}\Box B_{(n_2,\ldots,n_k)},
       \partial(J_{2r}\Box B_{(n_2,\ldots,n_k)})\bigr)\to(G,*),\\
  g&:\bigl(J_{2s}\Box B_{(n_2,\ldots,n_k)},
       \partial(J_{2s}\Box B_{(n_2,\ldots,n_k)})\bigr)\to(G,*).
 \end{aligned}
\]
For $k=1$, use the one-vertex box in place of $B_{(n_2,\ldots,n_k)}$.
Let
\[
 \iota_f:J_{2r}\hookrightarrow J_{2(r+s)},
 \qquad
 \iota_g:J_{2s}\hookrightarrow J_{2(r+s)}
\]
be the first-coordinate inclusions given by $\iota_f(t)=t$ and
$\iota_g(t)=2r+t$, respectively.
For the signed line chains,
\[
 c_{J_{2(r+s)}}=(\iota_f)_\#c_{J_{2r}}+(\iota_g)_\#c_{J_{2s}}.
\]
By bilinearity of the cross product,
\begin{equation}\label{eq:top-concat-chain}
 \begin{aligned}
 C_{(r+s,n_2,\ldots,n_k)}
 &=(\iota_f\Box\id_{B_{(n_2,\ldots,n_k)}})_\#C_{(r,n_2,\ldots,n_k)}\\
 &\quad +(\iota_g\Box\id_{B_{(n_2,\ldots,n_k)}})_\#C_{(s,n_2,\ldots,n_k)}.
 \end{aligned}
\end{equation}
The restrictions of $f\#_1g$ along $\iota_f\Box\id$ and
$\iota_g\Box\id$ are $f$ and $g$, respectively, by
\eqref{eq:concat}. Therefore
\[
 \begin{aligned}
 (f\#_1g)_\#C_{(r+s,n_2,\ldots,n_k)}
 &=f_\#C_{(r,n_2,\ldots,n_k)}
   +g_\#C_{(s,n_2,\ldots,n_k)},\\
 h_k^C([f]\cdot[g])
 &=h_k^C([f])+h_k^C([g]).
 \end{aligned}
\]

Let $f:(G,*)\to(H,*)$ be a based digraph map. For a standard
representative $f_0:(B_{\mathbf n},\partial B_{\mathbf n})\to(G,*)$ we have
\[
 \begin{aligned}
 h_k^C\bigl(f_*[f_0]\bigr)
 &=\bigl[(f\circ f_0)_\#C_{\mathbf n}\bigr]\\
 &=\bigl[f_\#((f_0)_\#C_{\mathbf n})\bigr]
 =H_k(f)\bigl(h_k^C([f_0])\bigr).
 \end{aligned}
\]
This proves naturality.
\end{proof}

For a standard representative $f$ on $J_{2n}$, the chain
$f_\#c_{J_{2n}}$ is the signed sum of the arrows traversed by $f$
in increasing parameter order
\cite[Eq.~(4.12)]{GLMY}. For a finite weakly connected digraph
$(G,*)$, $h_1^C$ is therefore the composite of $\alpha_1^G$ from
\cref{prop:pi1} with the GLMY abelianization
$\pi_1(G,*)\to H_1(G;\Z)$ \cite[Theorem~4.23]{GLMY}.
For infinite weakly connected $G$, take the directed colimit over
finite weakly connected induced subdigraphs containing $*$.
Every finite vertex set lies in one of these subdigraphs, obtained
by adjoining paths to $*$ in the underlying undirected graph.
Loops, $C$-homotopies and path chains have finite support, and
abelianization commutes with directed colimits.

\section{Standard spheres and a six-vertex example}\label{sec:examples}

We use the cubical nerve to compute the sphere homotopy groups of
standard digraph spheres and of a six-vertex digraph. The Hurewicz
homomorphism identifies a generator in the sphere dimension.

\subsection{Comparison with the cubical nerve}

For a digraph $X$, its \emph{cubical nerve} $N_1X$ is given by
\[
 (N_1 X)_r=\Hom(I_1^{\Box r},X),
 \qquad I_1=[0\to1],\qquad r\ge0,
\]
where $I_1^{\Box 0}$ is the one-vertex digraph and $\Hom$ denotes
digraph maps. A coordinate map $\theta:I_1^{\Box m}\to I_1^{\Box r}$
takes an $r$-cube $\sigma$ to the $m$-cube $\sigma\circ\theta$.
A face map inserts a coordinate equal to $0$ or $1$. A degeneracy
map deletes one coordinate. A connection map replaces two adjacent
coordinates by their maximum or minimum. A coordinate permutation
reorders the inputs of $\sigma$
\cite[Section~2.1, Definition~7.7 and Lemma~7.10]{KishimotoTong}.

To form the realization, take one copy of $[0,1]^r$ for each
$r$-cube $\sigma$. Use the same coordinate formula for $\theta$
on real coordinates and impose the relation
\[
 (\sigma\circ\theta,u)\sim(\sigma,\theta(u)),
 \qquad u\in[0,1]^m.
\]
We write $|N_1X|$ for the quotient by the equivalence relation generated
by the displayed pairs for face and degeneracy maps, and $[\sigma,u]$
for the class of $(\sigma,u)$ \cite[Definition~2.11]{KishimotoTong}.
The realization $|N_1X|_\gamma$ of Eldridge et al. adds the displayed
pairs for connection maps
\cite[Definition~2.1.6]{EldridgeEtAl}. For example, for a $1$-cube
$\sigma$, the point $(s,t)$ in the square labeled by
$\sigma\circ\max$ is identified with $\max(s,t)$ in the interval
labeled by $\sigma$. The resulting quotient map
$|N_1X|\to|N_1X|_\gamma$ is a homotopy equivalence
\cite[Definition~2.16 and Proposition~2.17]{KishimotoTong}.

\begin{proposition}\label{prop:sphere-grid-realization}
For every based digraph $(X,*)$ and every $j\ge1$, there are
natural group isomorphisms
\begin{equation}\label{eq:standard-sphere-grid-realization}
 \pC_j(X,*)\cong\pi_j^\delta(N_1X,*)\cong\pi_j(|N_1X|,*).
\end{equation}
Here $\pi_j^\delta(N_1X,*)$ is the group of spherical $j$-grids in
$N_1X$ modulo grid homotopy.
\end{proposition}
\begin{proof}
By \cref{rem:sphere-condition-standardization}, all oriented boxes
of positive side length may be used for representatives and common sources.
For an arrow $e$ of the $i$th factor $I_{m_i}$,
let $\iota_e:I_1\to I_{m_i}$ send $0$ to its source and $1$ to its
target. The grid associated to
$f:(B_{\mathcal I},\partial B_{\mathcal I})\to(X,*)$ is
\[
 f\longmapsto
 \bigl(f\circ(\iota_{e_1}\Box\cdots\Box\iota_{e_j})\bigr)_{(e_1,\ldots,e_j)}.
\]
The tuple records the restrictions of $f$ to all products of one
edge from each factor. These are the elementary cubes of the grid
\cite[Lemma~3.3 and Definition~3.12]{KishimotoTong}. Conversely,
the restrictions to a common face of two elementary cubes agree,
so their vertex values define a unique map on $B_{\mathcal I}$.
Every box arrow lies in an elementary cube, and a spherical grid has
value $*$ at each boundary vertex. The glued map is therefore a map
of pairs. Both assignments commute with postcomposition by based
digraph maps.

Under this correspondence the generating moves of grid homotopy
\cite[Definitions~3.5, 3.8, 3.9 and~3.12 and Lemma~3.7]{KishimotoTong}
have the forms
\[
 \begin{aligned}
  0\text{-homotopy}:\quad&
  f\simeq_0 f\circ\rho,\quad
  \rho:B_{\mathcal K}\to B_{\mathcal I}\text{ a product shrinking},\\
  1\text{-homotopy}:\quad&
  H:(B_{\mathcal I}\Box I_s,\partial B_{\mathcal I}\Box I_s)
      \longrightarrow(X,*),\\[-2pt]
  &H(-,0)=f,\qquad H(-,s)=g .
 \end{aligned}
\]
For a $0$-move, collapse the edge between $t$ and $t+1$ in a
coordinate line when the slices at $t$ and $t+1$ agree at every
choice of the other coordinates. Conversely, a coordinate shrinking
collapses precisely the edges whose endpoints have equal images
and preserves the orientation of every remaining edge. It therefore
factors into such single-edge collapses. For a $1$-move, place the
time factor last. An arrow $t\to t'$ of $I_s$ gives the pointwise
comparison $H(-,t)\toeq H(-,t')$ on $B_{\mathcal I}$, and each
slice has value $*$ on $\partial B_{\mathcal I}$.

A one-step $C$-homotopy with common source $B_{\mathcal K}$ and
shrinkings $\rho$ and $\sigma$ as in \cref{def:direct-sphere-C} yields
\[
 f\simeq_0 f\circ\rho\simeq_1 g\circ\sigma\simeq_0 g.
\]
Conversely, \eqref{eq:extension-C} handles the $0$-moves, and
identity shrinkings give the $C$-homotopy steps for the relative
comparisons of the $1$-moves. Hence
$\pC_j(X,*)\to\pi_j^\delta(N_1X,*)$ is a natural bijection.

Extend $f,g$ to the same standard box $B_{\nn}$, with
$\nn=(n_1,\ldots,n_j)$, by \cref{prop:tail-normal-form}.
The grid product of \cite[Section~3, following Definition~3.12]{KishimotoTong}
is the following map on $B_{2\nn}$:
\[
 (f\cdot g)(x)=
 \begin{cases}
  f(x),&0\le x_i\le 2n_i\text{ for all }i,\\
  g(x-2\nn),&2n_i\le x_i\le 4n_i\text{ for all }i,\\
  *,&\text{otherwise}.
 \end{cases}
\]
The formulas agree on block faces since $f,g$ are based on the boundary.
Define $p_-,p_+:B_{(2n_1,\ldots,2n_{j-1},n_j)}\to B_{\nn}$ by
\[
 \begin{aligned}
 (p_-(x))_i&=\min(x_i,2n_i),&
 (p_+(x))_i&=\max(0,x_i-2n_i) &&(i<j),\\
 (p_-(x))_j&=x_j,&(p_+(x))_j&=x_j .
 \end{aligned}
\]
Both are product shrinkings since the shifts are even. The boundary
values give
\[
 f\cdot g=(f\circ p_-)\#_j(g\circ p_+),\qquad
 [f\circ p_-]=[f],\quad [g\circ p_+]=[g].
\]
Here the last two equalities follow from \eqref{eq:extension-C}.
For $j\ge2$, interchange with coordinate $1$ as in \cref{prop:abelian}
shows that $\#_j$ induces the group product. For $j=1$, both
shrinkings are identities and $\#_1$ is the defining product.
This proves the first isomorphism in
\eqref{eq:standard-sphere-grid-realization}. The second follows from
\cite[Theorem~1.1 and Lemma~7.10]{KishimotoTong}.
\end{proof}

\subsection{Standard digraph spheres}\label{subsec:standard-examples}

Eldridge et al. computed the cubical realizations of the boundary
of $J_4\Box J_4$ and its out-closure using boundary strips
\cite[Example~4.0.9]{EldridgeEtAl}. We use the coordinate faces of
the interior box.

We use one cone construction for the standard spheres and the
six-vertex example. Let $A$ be a nonempty induced subdigraph of $K$
such that every arrow in $K$ starting in $A$ ends in $A$. Adjoin the
inverse cone $C^+(A)$ with a new vertex $v$ and arrows $v\to x$ for
$x\in V_A$.
In a directed cube, every vertex is reachable from its initial vertex
$\mathbf 0=(0,\ldots,0)$. For a cube $\sigma$
in $K\cup_A C^+(A)$,
\[
 v\in\operatorname{im}\sigma
 \quad\Longrightarrow\quad \sigma(\mathbf 0)=v
 \quad\Longrightarrow\quad
 \operatorname{im}\sigma\subseteq V_{C^+(A)}.
\]
The first implication uses that no arrow enters $v$. The second
uses the condition on arrows starting in $A$. Hence
\[
 N_1(K\cup_A C^+(A))=N_1K\cup_{N_1A}N_1C^+(A).
\]

For a one-step direct homotopy $H:X\Box I_1\to Y$, the prism
homotopy is given by
\[
 \begin{gathered}
  \widehat H_\sigma(w,\varepsilon)=H(\sigma(w),\varepsilon),
    \quad \sigma\in(N_1X)_r,\quad \varepsilon\in\{0,1\},\\
  |N_1X|\times[0,1]\longrightarrow |N_1Y|,\qquad
    ([\sigma,u],t)\longmapsto[\widehat H_\sigma,(u,t)] .
 \end{gathered}
\]
The prism respects faces and degeneracies in
the spatial coordinates and has endpoint maps induced by
$H(-,0)$ and $H(-,1)$. The constant map at $v$ points to the
identity of $C^+(A)$, so the prism contracts $|N_1C^+(A)|$.
The inclusions of $|N_1A|$ into $|N_1K|$ and $|N_1C^+(A)|$ are
inclusions of CW subcomplexes \cite[Lemma~2.12]{KishimotoTong}.
The realized union is therefore a homotopy pushout. Thus
\begin{equation}\label{eq:standard-sphere-nerve-union}
 |N_1K|\simeq\{\mathrm{pt}\}
 \quad\Longrightarrow\quad
 |N_1(K\cup_A C^+(A))|\simeq\Sigma |N_1A|,
\end{equation}
where $\Sigma$ is the unreduced suspension.

\begin{proof}[Proof of \cref{cor:standard-sphere-separation}]
Use the decomposition of $\Sn$ in \eqref{eq:cone-decomposition}.
In each interior line factor, send one endpoint to its neighbor
and fix every other vertex. The fold and the identity are related
by one pointwise comparison, in the direction of the edge between
the endpoint and its neighbor. Iterating the folds in each factor
gives $|N_1K_{\nn}|\simeq\{\mathrm{pt}\}$. For
$a\in\{1,2n_i-1\}$, let $F_{i,a}$ be the induced face of $K_{\nn}$
with $i$th coordinate $a$. We have $A_{\nn}=\bigcup_{i,a}F_{i,a}$.
Since $a$ is a sink in $K_i$,
\[
 x\in V_{F_{i,a}},\quad x\to y\text{ in }K_{\nn}
 \quad\Longrightarrow\quad y\in V_{F_{i,a}}.
\]

For cubes $\sigma$ in $A_{\nn}$, reachability and the sink condition give
\[
 \sigma(\mathbf 0)\in V_{F_{i,a}}
 \quad\Longrightarrow\quad \operatorname{im}\sigma\subseteq V_{F_{i,a}},
 \qquad
 |N_1A_{\nn}|=\bigcup_{i,a}|N_1F_{i,a}|.
\]
For a collection of faces indexed by $(i_l,a_l)$,
\[
 \begin{aligned}
  \bigcap_l |N_1F_{i_l,a_l}|
   &=\left|N_1\left(\bigcap_l F_{i_l,a_l}\right)\right|,\\
  \bigcap_l F_{i_l,a_l}
   &\cong
   \begin{cases}
    \varnothing,&i_l=i_{l'}\ \text{and }a_l\ne a_{l'}
                    \text{ for some }l,l',\\
    \displaystyle\mathop{\Box}_{i\notin\{i_l\}}K_i,&\text{otherwise}.
   \end{cases}
 \end{aligned}
\]
The empty product is a point. The nonempty intersections are
products of interior lines, so their realizations are contractible
by endpoint folds. A collection of faces has nonempty intersection
exactly when it contains at most one face for each $i$. The
simplicial nerve of this cover is therefore the join $(S^0)^{*k}$
of $k$ pairs of points. Since every nonempty intersection is
contractible, the completed nerve is the face poset of this
simplicial nerve. The CW nerve theorem
\cite[Theorem~2.13]{RamrasNerve} applies in all degrees.
Whitehead's theorem gives
\[
 |N_1A_{\nn}|\simeq |(S^0)^{*k}|\cong S^{k-1}.
\]
By \eqref{eq:standard-sphere-nerve-union},
\[
 |N_1\Sn|\simeq\Sigma|N_1A_{\nn}|\simeq S^k.
\]
Rotate $S^k$ to send the image of $*$ to $s_0$, and apply
\eqref{eq:standard-sphere-grid-realization}. In degree $k$,
$h_k^C$ maps $[q_{\nn}]$ to the generator $[z_{\nn}]$ by
\cref{thm:top-map,lem:standard-generator}. It is therefore a
surjection between infinite cyclic groups and hence an isomorphism.
\end{proof}

For $n_1=\cdots=n_k=2$, the standard quotient has vertex set
\[
 V_{S^k_{4\times\cdots\times4}}=\{*\}\sqcup\{1,2,3\}^k.
\]
For $x,y\in\{1,2,3\}^k$, its arrows are
\[
 \begin{aligned}
  x\to y
  &\iff \exists\,i:\ x_i=2,\ y_i\in\{1,3\},\
                x_l=y_l\text{ for }l\ne i,\\
  *\to x
  &\iff \exists\,i:\ x_i\in\{1,3\}.
 \end{aligned}
\]
See \cref{fig:standard-sphere-4x4} for $k=2$.

\subsection{The inverse suspension of the alternating four-cycle}\label{sec:six}

Let $R$ be the alternating four-cycle with vertices $r_0,r_1,r_2,r_3$ and arrows
\begin{equation}\label{eq:rim}
 r_0\to r_1,
 \qquad r_2\to r_1,
 \qquad r_2\to r_3,
 \qquad r_0\to r_3.
\end{equation}
The inverse suspension of $R$, based at $a$, is given by
\[
 V_G=V_R\sqcup\{a,b\},\qquad
 E_G=E_R\sqcup\{a\to r,\ b\to r:r\in V_R\}.
\]

\begin{figure}[ht]
\centering
\begin{tikzpicture}[
>=Stealth,
x=1cm,y=1cm,
edge/.style={->,thin,shorten <=3pt,shorten >=3pt},
crossing edge/.style={edge,
preaction={draw=white,line width=2.2pt,-}}
]
\coordinate (a)  at (0, 2.65);
\coordinate (b)  at (0,-2.65);
\coordinate (r0) at (-2.25,0);
\coordinate (r1) at (-.45, .72);
\coordinate (r2) at (2.25,0);
\coordinate (r3) at (.45,-.72);

\draw[edge] (a)--(r0);
\draw[edge] (a)--(r1);
\draw[edge] (a)--(r2);
\draw[edge] (b)--(r0);
\draw[edge] (b)--(r2);
\draw[edge] (b)--(r3);
\draw[edge,dashed] (b)--(r1);

\draw[edge] (r0)--(r1);
\draw[edge] (r2)--(r1);
\draw[edge] (r2)--(r3);
\draw[edge] (r0)--(r3);

\draw[crossing edge] (a)--(r3);

\node[graph vertex] at (a) {};
\node[graph vertex] at (b) {};
\node[graph vertex] at (r0) {};
\node[graph vertex] at (r1) {};
\node[graph vertex] at (r2) {};
\node[graph vertex] at (r3) {};
\node[font=\small] at (.24, 2.86) {$a$};
\node[font=\small] at (.24,-2.86) {$b$};
\node[font=\small] at (-2.58,.18) {$r_0$};
\node[font=\small] at (-.20, .98) {$r_1$};
\node[font=\small] at (2.58,.18) {$r_2$};
\node[font=\small] at (.72,-.98) {$r_3$};
\end{tikzpicture}
\caption{The inverse suspension $G$ of the alternating four-cycle $R$.
The rim $R$ is drawn as the equator,
with the two sources $a,b$ at the upper and lower poles.  The basepoint is $a$.}
\label{fig:six}
\end{figure}
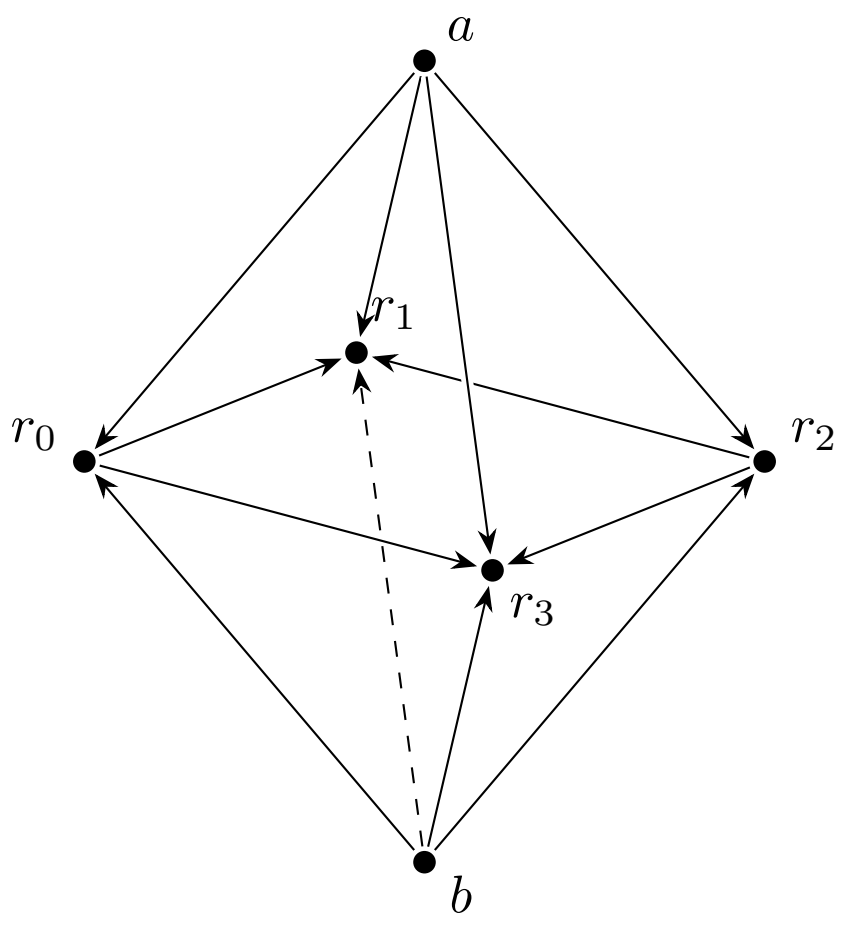

\begin{definition}\label{def:prefix-chain}
For $v\in\{a,b\}$ and $p\ge0$, define the \emph{prefix}
$v*u\in\mathcal R_{p+1}(G;\Z)$ of $u\in\mathcal R_p(R;\Z)$ by
setting $v*e_{x_0\cdots x_p}=e_{v x_0\cdots x_p}$ on elementary
paths and extending linearly. For $p\ge1$, the boundary formula in
the regular path module gives
\begin{equation}\label{eq:prefix-boundary-degree-two}
 \partial(v*u)=u-v*(\partial u).
\end{equation}
\end{definition}

Define
\begin{equation}\label{eq:rim-chain}
 c_R=e_{r_0r_1}-e_{r_2r_1}+e_{r_2r_3}-e_{r_0r_3},
 \qquad
 w=b*c_R-a*c_R.
\end{equation}
For $u\in\Omega_1(R;\Z)$ with coefficients $(u_0,u_1,u_2,u_3)$
in the order of \eqref{eq:rim},
\[
 \partial u=0
 \iff u_0+u_1=u_1+u_2=u_2+u_3=u_3+u_0=0
 \iff u\in\Z c_R .
\]
The rim has no allowed two-paths. In $G$, every allowed two-path
starts at a pole and follows a rim edge, and its middle face is
allowed because both poles point to every rim vertex. Thus
\[
 \begin{aligned}
  \Omega_2(R;\Z)&=0,\qquad H_1(R;\Z)=\Z[c_R],\\
  \Omega_2(G;\Z)&=\mathcal A_2(G;\Z)
    =a*\Omega_1(R;\Z)\oplus b*\Omega_1(R;\Z),\qquad
    \Omega_3(G;\Z)=0 .
 \end{aligned}
\]
For $x,y\in\Omega_1(R;\Z)$, apply \eqref{eq:prefix-boundary-degree-two}.
The three terms $x+y$, $a*(\partial x)$ and $b*(\partial y)$
have disjoint arrow supports, so
\[
 \begin{aligned}
  \partial(a*x+b*y)&=x+y-a*(\partial x)-b*(\partial y),\\
  \partial(a*x+b*y)=0
   &\iff x+y=0,\quad \partial x=\partial y=0\\
   &\iff a*x+b*y\in\Z w.
 \end{aligned}
\]
Consequently,
\begin{equation}\label{eq:oct-H2}
 H_2(G;\Z)=\Z[w].
\end{equation}

\begin{proposition}\label{prop:octahedral-sphere-class}
There is a digraph map of pairs
\[
 F:\bigl(J_4\Box J_4,\partial(J_4\Box J_4)\bigr)\longrightarrow(G,a)
\]
such that
\begin{equation}\label{eq:F-image}
 F_\#(c_{J_4}\times c_{J_4})=w.
\end{equation}
By \cref{thm:top-map} and \eqref{eq:oct-H2},
\[
 h_2^C([F])=[w]\ne0.
\]
\end{proposition}

\begin{proof}
Define $F(i,j)$ by
\begin{equation}\label{eq:F-table}
\begin{array}{c|ccccc}
 i\backslash j&0&1&2&3&4\\\hline
 0&a&a&a&a&a\\
 1&a&r_1&r_1&r_1&a\\
 2&a&r_2&b&r_0&a\\
 3&a&r_3&r_3&r_3&a\\
 4&a&a&a&a&a
\end{array}
\end{equation}
By \eqref{eq:F-table}, the boundary maps to $a$. Arrows between
the boundary and interior point inward, and $(2,2)$ is a source.
For a box arrow $u\to v$ with $F(u)\ne F(v)$, the table gives
\[
 \begin{cases}
 F(u)=a,\quad F(v)\in V_R, & u\text{ on the boundary},\\
 F(u)=b,\quad F(v)\in V_R, & u=(2,2),\\
 F(u)\in\{r_0,r_2\},\quad F(v)\in\{r_1,r_3\}, & \text{otherwise}.
 \end{cases}
\]
All three cases give arrows in $G$, so $F$ is a digraph map of pairs.

An allowed two-path in $G$ has middle vertex $r_2$ or $r_0$.
Their only preimages under $F$ are $(2,1)$ and $(2,3)$,
respectively. Thus
the shuffle terms with nonzero image are precisely
\[
 (2,j)\to(2,u)\to(v,u),
 \qquad (u,j)\in\{(1,0),(1,2),(3,2),(3,4)\},\quad v\in\{1,3\}.
\]
In each term the second coordinate changes first. By
\eqref{eq:two-shuffles}, its coefficient is minus the product of
the two edge coefficients in $c_{J_4}=e_{01}-e_{21}+e_{23}-e_{43}$.
Grouping the images by $u=1$ and $u=3$, respectively, gives
\begin{align*}
 F_\#(c_{J_4}\times c_{J_4})
 &= (e_{ar_2r_1}-e_{br_2r_1})-(e_{ar_2r_3}-e_{br_2r_3})\\
 &\quad +(e_{br_0r_1}-e_{ar_0r_1})-(e_{br_0r_3}-e_{ar_0r_3})\\
 &=b*c_R-a*c_R=w. \qedhere
\end{align*}
\end{proof}

\begin{corollary}\label{cor:six-vertex-separation}
For the inverse suspension $(G,a)$ of \cref{fig:six} and any
basepoint $s_0\in S^2$,
\[
 \pC_j(G,a)\cong\pi_j(S^2,s_0)\qquad(j\ge1).
\]
In particular, $\pC_1(G,a)=\pi_1(G,a)=0$. For $F$ as in
\cref{prop:octahedral-sphere-class}, the Hurewicz homomorphism is
\[
 h_2^C:\pC_2(G,a)=\Z[F]\xrightarrow{\ \cong\ } H_2(G;\Z)=\Z[w],
 \qquad
 [F]\longmapsto[w].
\]
The original GLMY homotopy groups satisfy
$\pi_j^{\mathrm{GLMY}}(G,a)=0$ for every $j\ge2$.
\end{corollary}
\begin{proof}
For a cube $\sigma:I_1^{\Box r}\to R$, each vertex of $I_1^{\Box r}$
lies on a directed path from $\mathbf 0$ to $\mathbf 1=(1,\ldots,1)$.
As $R$ has no allowed
two-paths,
\[
 \operatorname{im}\sigma\subseteq\{\sigma(\mathbf 0),\sigma(\mathbf 1)\},
 \qquad \sigma(\mathbf 0)\toeq\sigma(\mathbf 1).
\]
Hence $N_1R=\bigcup_{r\to s\in E_R}N_1(r\to s)$, where
$(r\to s)$ denotes the corresponding edge subdigraph. The realized
edge nerves are contractible CW subcomplexes. Adjacent ones meet at
their common endpoint, and opposite ones are disjoint. Thus every
nonempty intersection is contractible and the simplicial nerve is a
four-cycle. The same CW nerve theorem gives $|N_1R|\simeq S^1$.
In the cone construction of \eqref{eq:standard-sphere-nerve-union},
take $A=R$, $K=G\setminus\{b\}$ and $v=b$. The arrows starting
in $R$ stay in $R$, and $a\toeq x$ for all $x\in V_K$ gives
$|N_1K|\simeq\{\mathrm{pt}\}$ by the same prism used for the cone. Hence
\[
 |N_1G|\simeq\Sigma |N_1R|\simeq S^2.
\]
Rotate $S^2$ to send the image of $a$ to $s_0$.
Equation~\eqref{eq:standard-sphere-grid-realization} gives the
sphere homotopy groups, and \cref{prop:pi1} identifies
$\pC_1(G,a)$ with $\pi_1(G,a)$. In degree two, $h_2^C$ is a map between
infinite cyclic groups and is surjective by
\cref{prop:octahedral-sphere-class,eq:oct-H2}, so it is an
isomorphism. The original GLMY homotopy groups in degrees $j\ge2$
vanish by \cref{thm:GLMY-higher-vanishing}.
\end{proof}

\section{A loop digraph model for sphere homotopy groups}\label{sec:tail}

The six-vertex example distinguishes sphere homotopy groups from the
original GLMY homotopy groups. We now recover the sphere homotopy groups
by iterating a loop digraph whose arrows compare loops pointwise.

Li, Muranov and Yau compare the original loop digraph $L(G,*)$
with the replenished loop digraph by inclusion, and the replenished
loop digraph with the simple loop digraph by a based deformation
retraction \cite[Section~4 and Appendix~A]{LiMuranovYauSpheres}.
The inclusion and the retraction induce natural bijections on components
and natural isomorphisms on the positive-degree GLMY homotopy groups.
In particular, \cref{thm:GLMY-higher-vanishing} gives
$\pi_1(L(G,*),e_{G,*})=0$, and the other two loop digraphs also
have trivial fundamental groups at their distinguished basepoints.
By \cref{prop:pi1}, each of the three loop digraphs also has trivial
first sphere homotopy group at its distinguished basepoint.
A recursion identifying $\pC_2(G,*)$ with the first sphere homotopy
group of any of these loop digraphs for every $(G,*)$ would force
$\pC_2(G,*)=0$, contrary to \cref{cor:six-vertex-separation}.

\begin{definition}\label{def:J-infinity}
The \emph{infinite alternating line} $J_\infty$ has vertex set
$\{0,1,2,\ldots\}$ and, for each $r\ge0$, arrows
\[
 2r\to 2r+1,
 \qquad
 2r+2\to2r+1.
\]
These are all the arrows of $J_\infty$.
\end{definition}

\begin{definition}\label{def:tail-loop}
Let $(G,*)$ be a based digraph.  The \emph{pointwise-tail loop digraph}
$\Lt(G,*)$ has the following vertices, arrows and basepoint.
\begin{enumerate}
\item A vertex is a based digraph map
\[
 \gamma:(J_\infty,0)\longrightarrow(G,*)
\]
that is equal to the basepoint for all sufficiently large $t$.
\item Two distinct vertices $\gamma$ and $\delta$ satisfy
$\gamma\to\delta$ if and only if
\[
 \gamma(t)\toeq\delta(t)
 \qquad\text{for every }t\ge0.
\]
\item The basepoint is the constant sequence $t\mapsto *$.
\end{enumerate}
\end{definition}

A based digraph map $f:(G,*)\to(H,*)$ induces the map
\[
 \Lt(f):\Lt(G,*)\longrightarrow\Lt(H,*),
 \qquad \gamma\longmapsto f\circ\gamma.
\]
Thus $\Lt$ is a functor on based digraphs.

\begin{definition}\label{def:finite-loop-stage}
For a based digraph $(G,*)$ and $N\ge2$, let $\mathcal L_N(G,*)$
have the maps of pairs $(J_{2N},\partial J_{2N})\to(G,*)$ as vertices.
The basepoint is the constant loop. For distinct vertices $\gamma$
and $\delta$, there is an arrow $\gamma\to\delta$ if and only if
$\gamma(t)\toeq\delta(t)$ for all $0\le t\le2N$.
Extension by $*$ gives the map
\[
 \jmath_N:\mathcal L_N(G,*)\hookrightarrow\Lt(G,*),\qquad
 \jmath_N(\gamma)(t)=
 \begin{cases}
   \gamma(t),&0\le t\le 2N,\\
   *,&t>2N.
 \end{cases}
\]
The map $\jmath_N$ identifies $\mathcal L_N(G,*)$ with the induced
based subdigraph on the loops that equal $*$ for every $t\ge2N$.
\end{definition}

For $M\ge N$, the terminal shrinking $p_{M,N}$ gives an inclusion
$\mathcal L_N(G,*)\hookrightarrow\mathcal L_M(G,*)$ by precomposition,
since $\jmath_M(\gamma\circ p_{M,N})=\jmath_N(\gamma)$. Every loop lies
in one of these induced subdigraphs, so
\begin{equation}\label{eq:tail-colimit}
 \Lt(G,*)=\underset{N\ge2}{\operatorname{colim}}\,\mathcal L_N(G,*).
\end{equation}

Let $\mathbf S$ be the subcategory of shrinkings given by the
right-shrinking sequence in \cite[Example~1.2.12(4)]{EldridgeEtAl}.
Eldridge et al. define $\Omega_{\mathbf S}(G)$ as a colimit of
based mapping digraphs \cite[Definition~1.4.1]{EldridgeEtAl}.
Restrict to the cofinal even lengths $2N$ as in
\cref{rem:cubical-comparison}. A based map
$(J_{2N}/\partial J_{2N},*)\to(G,*)$ is equivalent to a map of pairs
$(J_{2N},\partial J_{2N})\to(G,*)$. Removing the reflexive arrows
leaves the pointwise comparisons between distinct maps, as in
$\mathcal L_N(G,*)$. The transitions are given by precomposition
with $p_{M,N}$. Hence \eqref{eq:tail-colimit} identifies $\Lt(G,*)$
naturally as the based digraph obtained from $\Omega_{\mathbf S}(G)$
by deleting its reflexive arrows.
The recursion in \cite[Proposition~1.4.5 and Corollary~1.4.6]{EldridgeEtAl}
uses such colimits of loop digraphs. We give a direct proof for
sphere representatives using slicing and evaluation.

\begin{lemma}\label{lem:finite-evaluation}
Let $(G,*)$ be a based digraph and $K$ a digraph with a subdigraph $A$.
For every $N\ge2$, evaluation gives a bijection
\begin{equation}\label{eq:exponential}
\begin{split}
 \ev_N:\;&\Hom((K,A),\mathcal L_N(G,*))\\
 &\longrightarrow
 \Hom\bigl((K\Box J_{2N},
 (A\Box J_{2N})\cup(K\Box\partial J_{2N})),(G,*)\bigr),
\end{split}
\end{equation}
with
\[
 \ev_N(\beta)(x,t)=\beta(x)(t),\qquad
 \ev_N^{-1}(f)(x)(t)=f(x,t).
\]
If $K$ is finite, the inclusions $\jmath_N$ induce a natural bijection
on homotopy classes relative to $A$,
\begin{equation}\label{eq:finite-evaluation-colimit}
 \underset{N\ge2}{\operatorname{colim}}
 \bigl[(K,A),\mathcal L_N(G,*)\bigr]
 \xrightarrow{\cong}
 \bigl[(K,A),\Lt(G,*)\bigr].
\end{equation}
\end{lemma}
\begin{proof}
The map $(x,t)\mapsto\beta(x)(t)$ preserves or collapses the arrows
in the $K$ coordinate precisely when
$\beta(x)(t)\toeq\beta(y)(t)$ for every $x\to y$ and $t$.
It preserves or collapses the arrows in the $J_{2N}$ coordinate
precisely when each $\beta(x)$ is a digraph map.
The boundary conditions require $\beta(a)$ to be the constant loop
for each $a\in V_A$ and $\beta(x)(0)=\beta(x)(2N)=*$ for each $x\in V_K$.
These conditions give the evaluation bijection.
If $K$ is finite, there is a common cutoff $2N$ for the finitely many
loops in the image of a map of pairs $(K,A)\to\Lt(G,*)$. The map
factors through the induced stage $\mathcal L_N(G,*)$, proving
surjectivity in \eqref{eq:finite-evaluation-colimit}.
For injectivity, suppose representatives from stages $N$ and $M$
are homotopic relative to $A$ in $\Lt(G,*)$. The homotopy has finite
domain $K\Box I_r$. Its loops have a common cutoff $2P$ with
$P\ge N,M$, so the two representatives become homotopic relative to
$A$ in $\mathcal L_P(G,*)$ after the stage inclusions.
\end{proof}

\begin{theorem}\label{thm:tail-recursion}
For every based digraph $(G,*)$ and $k\ge2$, there is a natural group
isomorphism
\[
 \Gamma_k:\pC_k(G,*)\xrightarrow{\cong}\pC_{k-1}(\Lt(G,*)).
\]
For a standard representative
\[
 f_0:(B_{\mathbf n'}\Box J_{2n_k},\partial(B_{\mathbf n'}\Box J_{2n_k}))
 \to(G,*),
\]
with $\mathbf n'=(n_1,\ldots,n_{k-1})$, its image is
$\Gamma_k([f_0])=[\beta]$, where
\[
 \beta:(B_{\mathbf n'},\partial B_{\mathbf n'})\longrightarrow\Lt(G,*),
 \qquad
 \beta(x)(t)=
 \begin{cases}
 f_0(x,t),&0\le t\le2n_k,\\
 *,&t>2n_k.
 \end{cases}
\]
\end{theorem}
\begin{proof}
Apply \cref{lem:finite-evaluation} with $(K,A)=(B_{\mathbf n'},\partial B_{\mathbf n'})$.
Evaluation preserves pointwise comparisons in both directions, so it
induces a bijection of relative homotopy classes at each finite stage.
Let $\beta_N:(B_{\mathbf n'},\partial B_{\mathbf n'})\to\mathcal L_N(G,*)$
be a map of pairs. For $\mathbf N'\ge\mathbf n'$ and $M\ge N$,
evaluation commutes with
terminal extensions in the first $k-1$ coordinates and in the last
coordinate, respectively:
\begin{equation}\label{eq:Gamma-shrinking}
 \ev_N(\beta_N\circ p_{\mathbf N',\mathbf n'})
 =\ev_N(\beta_N)\circ(p_{\mathbf N',\mathbf n'}\Box\id_{J_{2N}}),
\end{equation}
\begin{equation}\label{eq:cutoff}
 \ev_M\bigl(x\mapsto\beta_N(x)\circ p_{M,N}\bigr)
 =\ev_N(\beta_N)\circ(\id_{B_{\mathbf n'}}\Box p_{M,N}).
\end{equation}
By \cref{prop:tail-normal-form}, the colimit over $\mathbf n'$ of
the relative homotopy classes into $\Lt(G,*)$ is
$\pC_{k-1}(\Lt(G,*))$. Replace those classes by the colimit over
the finite stages from \eqref{eq:finite-evaluation-colimit} and
apply evaluation. The transition formulas
\eqref{eq:Gamma-shrinking} and \eqref{eq:cutoff} then give
\[
 \pC_{k-1}(\Lt(G,*))\cong
 \underset{\mathbf n',\,N}{\operatorname{colim}}
 \bigl[(B_{(\mathbf n',N)},\partial B_{(\mathbf n',N)}),(G,*)\bigr]
 \cong \pC_k(G,*).
\]
From right to left, the bijection $\Gamma_k$ takes $[f_0]$ to the class of
$\jmath_{n_k}\circ\ev_{n_k}^{-1}(f_0)$, as stated.

Use common terminal extensions for $f_0,f_1$ in coordinates $2,\ldots,k$,
and write $J_{2n_k}$ for their last factor. Slicing commutes with
concatenation in the first coordinate by the formula
\[
 \jmath_{n_k}\circ\ev_{n_k}^{-1}(f_0\#_1 f_1)=
 (\jmath_{n_k}\circ\ev_{n_k}^{-1}(f_0))\#_1(\jmath_{n_k}\circ\ev_{n_k}^{-1}(f_1)).
\]
Hence $\Gamma_k$ is a group isomorphism. For a based digraph map
$f:(G,*)\to(H,*)$, the slice of $f\circ f_0$ at $x$ is
$f\circ\beta(x)$, including on its constant tail since $f(*)=*$.
This gives the commutative square
\[
\begin{tikzcd}[column sep=large,row sep=large]
 \pC_k(G,*)\arrow[r,"\Gamma_k","\cong"']\arrow[d,"f_*"']
 &\pC_{k-1}(\Lt(G,*))\arrow[d,"(\Lt(f))_*"]\\
 \pC_k(H,*)\arrow[r,"\Gamma_k","\cong"']
 &\pC_{k-1}(\Lt(H,*)).
\end{tikzcd}
\qedhere
\]
\end{proof}

\begin{proposition}\label{prop:tail-components}
For every based digraph $(G,*)$, concatenation induces a group structure
on $\pi_0(\Lt(G,*))$, and extension by $*$ gives a natural group isomorphism
\[
 \pC_1(G,*)\cong\pi_0(\Lt(G,*)).
\]
\end{proposition}
\begin{proof}
For $K$ a single vertex and $A=\varnothing$, homotopy classes of
maps from $(K,A)$ to a based digraph are its components. Thus
\cref{lem:finite-evaluation,prop:tail-normal-form} give
\[
 \pi_0(\Lt(G,*))\cong
 \underset{N\ge2}{\operatorname{colim}}
 \bigl[(J_{2N},\partial J_{2N}),(G,*)\bigr]\cong\pC_1(G,*).
\]
The resulting bijection from $\pC_1(G,*)$ to $\pi_0(\Lt(G,*))$
sends a standard class $[f_0]$ on $J_{2m}$ to $[\jmath_m(f_0)]$.
For a second standard representative $f_1$ on $J_{2n}$,
\cref{prop:concat,thm:sphere-groups} show that
\[
 [\jmath_m(f_0)]\cdot[\jmath_n(f_1)]
 :=[\jmath_{m+n}(f_0\#_1f_1)]
\]
defines a group structure on $\pi_0(\Lt(G,*))$ for which this bijection is
an isomorphism.  Naturality follows, for a based digraph map
$f:(G,*)\to(H,*)$, from
\[
 \Lt(f)(\jmath_m(f_0))=\jmath_m(f\circ f_0). \qedhere
\]
\end{proof}

Define the iterates as based digraphs by
\[
 \Lt^0(G,*)=(G,*),\qquad
 \Lt^{r+1}(G,*)=\Lt\bigl(\Lt^r(G,*)\bigr)
 \quad(r\in\mathbb Z_{\ge0}).
\]
The basepoint of each positive iterate is the constant sequence on
$J_\infty$ at the preceding basepoint.

\begin{corollary}\label{cor:iterated-tail}
For every based digraph $(G,*)$ and every $k\ge1$, there are natural
group isomorphisms
\[
 \pC_k(G,*)\cong\pi_1(\Lt^{\,k-1}(G,*))
 \cong\pi_0(\Lt^{\,k}(G,*)),
\]
where the multiplication on components is given by \cref{prop:tail-components}.
\end{corollary}
\begin{proof}
Iterate \cref{thm:tail-recursion} to obtain
\[
 \pC_k(G,*)\cong\pC_1(\Lt^{\,k-1}(G,*)).
\]
For $k=1$ this identification is the identity. Now apply
\cref{prop:pi1,prop:tail-components} to $\Lt^{\,k-1}(G,*)$.
Naturality follows from that of the individual isomorphisms.
\end{proof}

Return to the six-vertex inverse suspension $(G,a)$.
Since $\pi_1(G,a)=0$, \cref{prop:pi1,prop:tail-components} give
$|\pi_0(L(G,a))|=|\pi_0(\Lt(G,a))|=1$.
With the constant loops as basepoints, we have
\[
 \pi_1(L(G,a),e_{G,a})=0,\qquad
 \pi_1(\Lt(G,a))\cong\pC_2(G,a)=\Z[F],
\]
by \cref{thm:GLMY-higher-vanishing,cor:six-vertex-separation,cor:iterated-tail}.
The generator $\alpha_1^{\Lt(G,a)}(\Gamma_2([F]))$ of
$\pi_1(\Lt(G,a))$ is represented by the based loop
\[
 J_4\longrightarrow\Lt(G,a),\qquad
 i\longmapsto\bigl(F(i,0),F(i,1),F(i,2),F(i,3),F(i,4),a,a,\ldots\bigr),
\]
with $F(i,t)$ given by row $i$ of \eqref{eq:F-table}. The two
endpoint rows give the constant loop at $a$.

\appendix
\section{GLMY path homology of standard digraph spheres}\label{app:homology}

The standard sphere presentations in Section~\ref{sec:spheres}
rely on \cref{thm:admissibility}. To prove the theorem, we reduce
the GLMY path complex of a standard sphere relative to its basepoint
to a tensor product of relative line complexes over $\Z$.

For an induced subdigraph $Y\subseteq X$, write
\[
 \Omega_*(X,Y;\Z)=\Omega_*(X;\Z)/\Omega_*(Y;\Z),
 \qquad H_*(X,Y;\Z)=H_*(\Omega_*(X,Y;\Z)).
\]
Fix $k\ge1$ and $\nn=(n_1,\ldots,n_k)$ with $n_i\ge2$.
Let $K_i$ be the subdigraph of $J_{2n_i}$ induced by
$\{1,\ldots,2n_i-1\}$, and let $E_i$ be its induced subdigraph on
$\{1,2n_i-1\}$.  Set
\[
 K_{\nn}=K_1\Box\cdots\Box K_k,
 \qquad
 A_{\nn}=\bigcup_{i=1}^k
 K_1\Box\cdots\Box E_i\Box\cdots\Box K_k.
\]
Every vertex of $E_i$ is a sink in $K_i$. For an arrow
$x\to y$ in $K_{\nn}$ with $x_i\in V_{E_i}$, we have $y_i=x_i$.
Thus the arrow lies in the $i$th product in the union defining
$A_{\nn}$. In particular, every arrow starting at a vertex of
$A_{\nn}$ stays in $A_{\nn}$, which is an induced subdigraph
of $K_{\nn}$.

The quotient $\Sn$ has the vertices of $K_{\nn}$ and the collapsed
boundary vertex $*$.  No arrow enters $*$, because $0$ and $2n_i$
are sources in $J_{2n_i}$. An arrow from the boundary to an interior
vertex $x$ must change an endpoint coordinate $0$ or $2n_i$ to
$1$ or $2n_i-1$, respectively. Conversely, each such coordinate
of $x$ gives an arrow ending at $x$ by restoring the neighboring
endpoint in the source. Therefore $*\to x$ precisely for
$x\in V_{A_{\nn}}$, and
\begin{equation}\label{eq:cone-decomposition}
 \Sn=K_{\nn}\cup_{A_{\nn}}C^+(A_{\nn}),
\end{equation}
where $C^+(A_{\nn})$ is the inverse cone
\cite[Definition~2.11]{Axioms}, obtained by adjoining $*$ and the
arrows $*\to x$ for all $x\in V_{A_{\nn}}$.

\begin{lemma}\label{lem:cone-chain}
The inclusion $K_{\nn}\hookrightarrow\Sn$ induces an isomorphism
of integral chain complexes
\[
 \Omega_*(K_{\nn},A_{\nn};\Z)
 \xrightarrow{\cong}
 \Omega_*(\Sn,C^+(A_{\nn});\Z).
\]
\end{lemma}
\begin{proof}
Every allowed path of positive degree in $\Sn$ either lies in
$K_{\nn}$ or starts at $*$ and subsequently lies in $A_{\nn}$.
Indeed, no arrow enters $*$, and every arrow starting in $A_{\nn}$
ends there. Split an allowed chain into the sum of its terms in
$K_{\nn}$ and its terms starting at $*$. The nonallowed regular
faces from the first chain avoid $*$. Those from the second
contain $*$, since deleting the initial $*$ leaves an allowed
path in $A_{\nn}$. The two sets of nonallowed faces have disjoint
supports. Thus the sum of the two chains is $\partial$-invariant
precisely when each chain is. The vertex basis gives the same
decomposition in degree zero. Therefore
\begin{equation}\label{eq:Omega-sum}
 \Omega_p(\Sn;\Z)
 =\Omega_p(K_{\nn};\Z)+\Omega_p(C^+(A_{\nn});\Z)
 \qquad(p\ge0).
\end{equation}
$K_{\nn}$ and $C^+(A_{\nn})$ are induced subdigraphs of $\Sn$
with intersection $A_{\nn}$. Their path-chain complexes satisfy
\[
 \Omega_p(K_{\nn};\Z)\cap\Omega_p(C^+(A_{\nn});\Z)
 =\Omega_p(A_{\nn};\Z).
\]
Thus the composite
\[
 \Omega_p(K_{\nn};\Z)\hookrightarrow\Omega_p(\Sn;\Z)
 \twoheadrightarrow\Omega_p(\Sn,C^+(A_{\nn});\Z)
\]
has kernel $\Omega_p(A_{\nn};\Z)$ and is surjective by
\eqref{eq:Omega-sum}. Both maps commute with $\partial$.
\end{proof}

We next contract the chain complex of the inverse cone relative to $*$.
Define $s:\mathcal R_p(C^+(A_{\nn});\Z)\to
\mathcal R_{p+1}(C^+(A_{\nn});\Z)$ on the elementary regular paths by
\[
 s(e_{x_0\cdots x_p})=
 \begin{cases}
  e_{*x_0\cdots x_p},&x_0\in V_{A_{\nn}},\\
  0,&x_0=*.
 \end{cases}
\]
Extend $s$ linearly. An allowed path starting in $A_{\nn}$ stays
there, so prefixing it with $*$ gives an allowed path. The operator
$s$ therefore sends allowed chains to allowed chains. For $p\ge1$,
alternating deletion in the regular path module gives
\[
 \partial e_{*x_0\cdots x_p}
 =e_{x_0\cdots x_p}-s\partial e_{x_0\cdots x_p}\quad (x_0\ne *),
 \qquad s\partial e_{*x_1\cdots x_p}=e_{*x_1\cdots x_p}.
\]
Hence $\partial s+s\partial=\id$ in positive degree. In degree zero,
\[
 (\partial s+s\partial)e_x=e_x-e_*
 \quad(x\in V_{A_{\nn}}),
 \qquad (\partial s+s\partial)e_*=0.
\]
For $u\in\Omega_p(C^+(A_{\nn});\Z)$, both $u$ and $s(\partial u)$
are allowed, so these identities imply
$\partial s(u)\in\mathcal A_p(C^+(A_{\nn});\Z)$.
Thus $s$ acts on the $\partial$-invariant path complex and sends
$e_*$ to zero. It descends to
$\Omega_*(C^+(A_{\nn}),\{*\};\Z)$, on which
\[
 \partial s+s\partial=\id.
\]
Thus $s$ contracts the relative cone complex.
The projection in the short exact sequence
\[
\begin{aligned}
0&\longrightarrow\Omega_*(C^+(A_{\nn}),\{*\};\Z)
 \longrightarrow\Omega_*(\Sn,\{*\};\Z)\\
 &\longrightarrow\Omega_*(\Sn,C^+(A_{\nn});\Z)
 \longrightarrow0
\end{aligned}
\]
induces an isomorphism on homology because the left-hand complex
is contractible. Combining this isomorphism with
\cref{lem:cone-chain} gives
\begin{equation}\label{eq:relative-reduction}
 H_q(\Sn,\{*\};\Z)\cong H_q(K_{\nn},A_{\nn};\Z)
 \qquad(q\ge0).
\end{equation}
The underlying undirected graph of $K_{\nn}$ is connected, and
$A_{\nn}$ is nonempty.  By \eqref{eq:cone-decomposition}, $\Sn$ is
weakly connected.  Hence
$H_0(\Sn;\Z)\cong\Z$, and
$H_q(\Sn;\Z)\cong H_q(\Sn,\{*\};\Z)$ for $q>0$.

\begin{lemma}\label{lem:shuffle-basis}
Let $X_1,\ldots,X_k$ be finite digraphs in which every arrow ends at
a sink.  The iterated path cross product is an isomorphism of
integral chain complexes
\[
 \bigotimes_{i=1}^k\Omega_*(X_i;\Z)
 \xrightarrow{\cong}\Omega_*(X_1\Box\cdots\Box X_k;\Z),
 \qquad
 \xi_1\otimes\cdots\otimes\xi_k\longmapsto
 \xi_1\times\cdots\times\xi_k,
\]
where the tensor product has the usual graded differential
\[
 \partial(\xi_1\otimes\cdots\otimes\xi_k)
 =\sum_{j=1}^k(-1)^{\sum_{i<j}|\xi_i|}
 \xi_1\otimes\cdots\otimes\partial\xi_j\otimes\cdots\otimes\xi_k
\]
for homogeneous factors.
\end{lemma}
\begin{proof}
If one factor has no vertices, both sides are zero.  Otherwise,
$\Omega_p(X_i;\Z)=0$ for $p\ge2$, and $\Omega_p(X_i;\Z)$ is freely
generated by the allowed $p$-paths for $p=0,1$.
In the box product the $i$th coordinate of an allowed path changes
at most once, since every arrow of $X_i$ ends at a sink.
Thus there are no allowed paths of length greater than $k$.

Fix $1\le p\le k$. A basis tensor
$\xi_1\otimes\cdots\otimes\xi_k$ of degree $p$ fixes a vertex in each
degree-zero factor and an arrow in each degree-one factor. The
initial product vertex uses the fixed vertices and the sources
of the $p$ arrows. The $p!$ orders of those coordinate arrows give
distinct allowed $p$-paths in the box product.  Conversely, each allowed $p$-path determines its initial vertex,
its coordinate arrows and their order. The allowed paths therefore
split into disjoint families indexed by the degree-$p$ basis tensors.  In degrees zero
and one, the cross product identifies the vertex and arrow bases,
respectively.

Now let $2\le p\le k$, fix a degree-$p$ basis tensor
$\xi_1\otimes\cdots\otimes\xi_k$, and take
\[
 u\in\mathcal A_p(X_1\Box\cdots\Box X_k;\Z)
\]
supported on its family of $p!$ paths.  Write
$I=\{i_1<\cdots<i_p\}$ for the coordinates that change.  For a
permutation $\sigma$ of $I$, let $a_\sigma\in\Z$ be the coefficient
in $u$ of the path whose coordinate changes occur in the order
$\sigma(i_1),\ldots,\sigma(i_p)$.

Consider the elementary path $e_{x_0\cdots x_p}$ with this order.
Deleting either endpoint leaves an allowed face.  Deleting $x_j$,
where $1\le j<p$, produces the regular face
\[
 e_{x_0\cdots\widehat{x_j}\cdots x_p}.
\]
Its consecutive vertices $x_{j-1}$ and $x_{j+1}$ differ in two
coordinates, so this face is nonallowed.  To obtain this face from an allowed path, one must insert a vertex
between $x_{j-1}$ and $x_{j+1}$. Inserting it elsewhere leaves the
nonallowed step unchanged.
The two possible routes are
\[
\begin{aligned}
 x_{j-1}&\xrightarrow{\sigma(i_j)}x_j
          \xrightarrow{\sigma(i_{j+1})}x_{j+1},\\
 x_{j-1}&\xrightarrow{\sigma(i_{j+1})}x_j'
          \xrightarrow{\sigma(i_j)}x_{j+1},
\end{aligned}
\]
where each arrow is labeled by the coordinate it changes, and
$x_j'$ is obtained from $x_{j-1}$ by changing coordinate
$\sigma(i_{j+1})$ to its value in $x_{j+1}$.  These are the only
possibilities, since each box-product arrow changes exactly one
coordinate.  Write $\sigma'$ for the permutation obtained by interchanging the
$j$th and $(j+1)$st entries. Replacing $x_j$ by $x_j'$ makes exactly
this interchange.  Both paths contribute the displayed
face with boundary sign $(-1)^j$, so its coefficient in $\partial u$
is
\[
 (-1)^j(a_\sigma+a_{\sigma'}).
\]
An interior face retains both endpoint vertices. At each coordinate
the endpoints determine the fixed vertex or the unique arrow from
the initial value to the final value. They therefore determine the
basis tensor. Each nonallowed face occurs in just one family of
paths, and it arises by deleting an interior vertex. Consequently,
\[
\begin{aligned}
 u\in\Omega_p(X_1\Box\cdots\Box X_k;\Z)
 \quad\Longleftrightarrow\quad
 &a_{\sigma'}=-a_\sigma\\
 &\text{for every adjacent interchange}.
\end{aligned}
\]
Let $\id$ denote the identity permutation of $I$, and compute
permutation signs relative to the increasing order $i_1<\cdots<i_p$.
Adjacent interchanges generate all permutations.  By the shuffle
sign convention of \cref{def:staircase-product},
\[
 a_\sigma=\operatorname{sgn}(\sigma)a_{\id},
 \qquad
 u=a_{\id}(\xi_1\times\cdots\times\xi_k).
\]
Conversely, these coefficients cancel every nonallowed face.
Thus the $\partial$-invariant chains supported on this family form
the infinite cyclic group generated by
$\xi_1\times\cdots\times\xi_k$.  Since the families have disjoint
supports, these cross products form a $\Z$-basis of $\Omega_p$.
The Leibniz identity for the cross product gives the tensor
differential in the statement. Thus the cross product is an
isomorphism of integral chain complexes.
\end{proof}

Apply \cref{lem:shuffle-basis} to $K_1,\ldots,K_k$.
For a family of paths indexed by a basis tensor, the initial vertex
is the same in every order of the coordinate arrows. If its $i$th
coordinate is in $V_{E_i}$, that coordinate is fixed, because the
vertices of $E_i$ are sinks. Hence a family lies in $A_{\nn}$ precisely
when some factor of its basis tensor is a vertex of $E_i$ in
coordinate $i$. The families have disjoint supports, and
$A_{\nn}$ is induced in $K_{\nn}$. Thus $\Omega_*(A_{\nn};\Z)$
corresponds under the shuffle basis to
\begin{equation}\label{eq:boundary-tensor}
 \sum_{i=1}^k
 \Omega_*(K_1;\Z)\otimes\cdots\otimes\Omega_*(E_i;\Z)
 \otimes\cdots\otimes\Omega_*(K_k;\Z).
\end{equation}
In each factor, split off the two vertex basis elements spanning
$\Omega_*(E_i;\Z)$. The tensor of the quotient maps kills the
summands containing an $\Omega_*(E_i;\Z)$ factor and maps the
tensor of the complementary summands isomorphically onto the
target. Its kernel is \eqref{eq:boundary-tensor}, giving the chain
isomorphism
\begin{equation}\label{eq:relative-tensor}
 \Omega_*(K_{\nn},A_{\nn};\Z)
 \cong\bigotimes_{i=1}^k\Omega_*(K_i,E_i;\Z).
\end{equation}

Since every arrow of $K_i$ ends at a sink, the relative complex
$\Omega_*(K_i,E_i;\Z)$ has only its arrow and vertex groups, in
degrees one and zero. Orient each edge in order of increasing vertex
index, taking the negative of its arrow basis element when necessary.
Then the boundary of the signed edge from $t$ to $t+1$ is
$e_{t+1}-e_t$, with the two endpoint vertices zero in the relative
complex. Thus the sum of the signed edges from $1$ through $j$
has boundary $e_j$, for $2\le j\le 2n_i-2$. Sending each
remaining vertex basis element $e_j$ to this edge sum defines a
section of the degree-one boundary.
In the signed edge basis, a relative one-chain is a cycle precisely
when the coefficients agree on the two sides of each interior vertex.
Thus all coefficients of a cycle are equal. The signed traversal of
the whole interval therefore generates the kernel over $\Z$.

The section splits off the degree-one kernel $\Z$. The remaining
two-term complex has an isomorphism as differential. Tensoring
these splittings gives one copy of $\Z$ in degree $k$. For any
other summand, choose a two-term factor in position $j$ and let
$s_j$ be its inverse differential, extended by zero in degree one.
By the tensor differential in \cref{lem:shuffle-basis}, the map
\[
 \xi_1\otimes\cdots\otimes\xi_k\longmapsto
 (-1)^{\sum_{i<j}|\xi_i|}
 \xi_1\otimes\cdots\otimes s_j(\xi_j)\otimes\cdots\otimes\xi_k
\]
contracts that summand.
Hence
\[
 H_q(K_{\nn},A_{\nn};\Z)\cong
 \begin{cases}
  \Z,&q=k,\\
  0,&q\ne k.
 \end{cases}
\]
Together with \eqref{eq:relative-reduction} and the connectedness of
$\Sn$, this proves \cref{thm:admissibility}.

\section*{Acknowledgements}
Jingyan Li was supported in part by the National Natural Science Foundation of China
(General Program, Grant No.~82573048) and the Beijing Natural Science Foundation
(International Scientists Project, Grant No.~IS25081).
Jingyan Li and Jie Wu were supported in part by the research start-up fund
of the Beijing Institute of Mathematical Sciences and Applications (BIMSA).

\enlargethispage{2\baselineskip}
\end{document}